\documentclass[a4wide,11pt]{article}
\usepackage[OT1,T1]{fontenc}
\usepackage{a4wide}
\usepackage{amsmath,amssymb}
\usepackage{mathrsfs}
\usepackage[ps,all]{xy}
\SelectTips{cm}{10}
\usepackage[pdfborder={0 0 0},pdfpagemode=None,pdfstartview=FitH]{hyperref}

\renewcommand\labelitemi{$\ \bullet$}
\renewcommand\theenumi{(\roman{enumi})}
\renewcommand\labelenumi{\theenumi}

\renewcommand\enumerate{\list\labelenumi
  {\setlength\leftmargin{0pt}\setlength\labelwidth{0pt}
  \usecounter{enumi}\def\makelabel##1{\kern\labelsep{##1}}}}
\renewcommand\itemize{\list\labelitemi
  {\setlength\leftmargin{0pt}\setlength\labelwidth{0pt}
  \def\makelabel##1{\kern\labelsep{##1}}}}
\renewenvironment{description}{\list{}
  {\setlength\leftmargin{0pt}\setlength\labelwidth{0pt}
  \let\makelabel\descriptionlabel}}
  {\endlist}

\newtheorem{theorem}{Theorem}
\newtheorem{lemma}[theorem]{Lemma}
\newtheorem{proposition}[theorem]{Proposition}

\newenvironment{proof}{\trivlist
  \item[\hskip\labelsep{\itshape Proof.}]\upshape}{\nobreak\noindent
  $\square$\endtrivlist}
\newenvironment{other}[1]{\refstepcounter{theorem}\trivlist
  \item[\hskip\labelsep{\itshape #1~\arabic{theorem}.}]\upshape}
  {\endtrivlist\bigbreak}
\newenvironment{other*}[1]{\trivlist
  \item[\hskip\labelsep{\itshape #1.}]\upshape}
  {\endtrivlist\bigbreak}

\DeclareMathOperator\Hom{Hom}
\DeclareMathOperator\Aut{Aut}
\DeclareMathOperator\Ext{Ext}
\DeclareMathOperator\im{im}
\DeclareMathOperator\coker{coker}
\DeclareMathOperator\id{id}
\DeclareMathOperator\hd{hd}
\DeclareMathOperator\soc{soc}
\DeclareMathOperator\Rep{Rep}
\DeclareMathOperator\wt{wt}
\DeclareMathOperator\Pol{Pol}
\DeclareMathOperator\Irr{Irr}
\DeclareMathOperator\spn{span}
\newcommand\dimvec{{\,\underline\dim\,}}
\newcommand\bsm{\begin{smallmatrix}}
\newcommand\esm{\end{smallmatrix}}
\newcommand\out{\text{out}}
\newcommand\iin{\text{in}}

\newcommand\MV{\mathcal{MV}}
\newcommand\IC{\mathrm{IC}}

\begin{document}
\title{Preprojective algebras and MV polytopes}
\author{Pierre Baumann and Joel Kamnitzer}
\date{}
\maketitle

\begin{abstract}
The purpose of this paper is to apply the theory of MV
polytopes to the study of components of Lusztig's nilpotent varieties.
Along the way, we introduce reflection functors for modules over
the non-deformed preprojective algebra of a quiver.
\end{abstract}

\section{Introduction}
\label{se:Intro}
Let $\mathfrak g$ be a simply-laced semisimple finite dimensional complex
Lie algebra. Part of the representation theory of $\mathfrak g$ is
described by Kashiwara's combinatorics of crystals \cite{Kashiwara95}. 
On the algebraic
side, crystals are implemented by special bases in representations, for
instance Lusztig's canonical and semicanonical bases 
\cite{Lusztig90a, Lusztig91, Lusztig00}. On the combinatorial
side, a crystal can be realized by several models, each of which provides a
concrete access to it; as an example, Anderson \cite{Anderson03} and the 
second author's \cite{Kamnitzer05,Kamnitzer07}
MV (Mirkovi\'c, Vilonen) polytopes form a model for the crystal
$B(-\infty)$ of the positive part $U(\mathfrak n)$ of $U(\mathfrak g)$.

The purpose of this paper is to apply the theory of MV polytopes to the
study of the semicanonical basis. This basis arises through a geometric
realisation of $U(\mathfrak n)$ by means of constructible functions on
certain varieties $\Lambda(\nu)$, called Lusztig's nilpotent varieties. 
Concretely, these varieties parametrize representations of the 
preprojective algebra constructed using the Dynkin graph of $\mathfrak g$. 
We are thus led to investigate the
relations between these varieties $\Lambda(\nu)$ and MV polytopes.

Our method is the following. For each chamber weight $\gamma$, we
construct a module $N(\gamma)$ over the preprojective algebra.  Then we 
define a constructible function $D_\gamma$ on the nilpotent
varieties whose value at a point $ T \in \Lambda(\nu) $ is 
$\dim \Hom(N(\gamma),T)$. On an irreducible component $Z$ of a nilpotent 
variety, the functions $D_\gamma$ admit generic values. We prove that these
generic values form the BZ (Berenstein, Zelevinsky) datum of an MV
polytope, say $\Pol(Z)$. The resulting bijection $Z\mapsto\Pol(Z)$ defines 
a labelling of the semicanonical basis by MV polytopes, and hence by
the crystal~$B(-\infty)$.

This indexation of the semicanonical basis by $B(-\infty)$ is not new;
indeed it was proved by Kashiwara and 
Saito \cite{Kashiwara-Saito97} (and we use their
result). However an important property here is the fact that the MV
polytope of an element $b\in B(-\infty)$ geometrically packs together
all the Lusztig parametrizations of $b$. Therefore our result can be
seen as the combinatorial facet of the relation between the
semicanonical basis and the PBW (Poincar\'e, Birkhoff, Witt) bases
of $U(\mathfrak n)$. A remarkable fact here is that this relation holds
even when the reduced decomposition used to define the PBW basis is not
adapted to any orientation of the graph of $\mathfrak g$.

To construct the modules $N(\gamma)$ and study the functions
$D_\gamma$, a key tool is reflection functors for modules over the
preprojective algebra. These functors extend the usual BGP
(Bernstein, Gelfand, Ponomarev) reflection functors
\cite{Bernstein-Gelfand-Ponomarev73} to the
preprojective framework; they are different from Crawley-Boevey
and Holland's reflection functors \cite{Crawley-Boevey-Holland98} for 
the deformed preprojective
algebra. Our initial motivation to define them was to understand
the meaning of Kashiwara and Saito's crystal reflections in terms
of the semicanonical basis.

The combinatorics of MV polytopes was originally developed by
Anderson and the second author, in order to describe the Mirkovi\'c-Vilonen
cycles in affine Grassmannians.  Affine Grassmannians and quiver varieties
are two geometric constructions of representations of $ \mathfrak g$ and 
they each give their own basis for representations.  One motivation
for the current paper was to relate the geometry of affine 
Grassmannians and nilpotent varieties, as part of a bigger project which 
seeks to understand and compare these bases.  In particular, in this paper
we have achieved a labelling of the irreducible components of nilpotent
varieties by MV polytopes and hence a natural bijection between these
components and MV cycles. Moreover, we have done so by constructing 
functions $ D_\gamma $ which are the direct analogs of functions introduced
by the second author in \cite{Kamnitzer05} for the study of MV cycles.  
In fact, the functions $D_\gamma$
define a stratification of the varieties $\Lambda(\nu)$. The strata
are indexed by the same pseudo-Weyl polytopes that index the GGMS
(Gelfand, Goresky, MacPherson, Serganova) strata of the affine
Grassmannian.

Another motivation of the current project is to develop a theory of 
MV polytopes for affine Kac-Moody Lie algebras.  From the MV cycles
perspective this is difficult, since the affine Grassmannian for Kac-Moody
groups is not yet available.  However, the theory of quivers extends nicely
to the Kac-Moody setting.  In a future work, we plan to extend the results
of this paper to quivers for affine Kac-Moody Lie algebras and obtain a
notion of MV polytopes in this setting.

In a related work \cite{Kamnitzer-Sadanand10}, the second author and 
C.~Sadanand analyzed explicity 
the case of $ \mathfrak g = \mathfrak{sl}_n $ and compared the results
of the current paper 
with a previous description of the components of nilpotent varieties
obtained by A.~Savage \cite{Savage07}.

We heartily thank A.~Braverman, C.~Gei\ss, B.~Keller, B.~Leclerc, 
C.~Sadanand, A.~Savage, J.~Schr\"oer, and P.~Tingley for
useful conversations. Gei\ss, Leclerc and Schr\"oer explained to us that
their work on the cluster properties of the semicanonical basis also
relies on the modules $N(\gamma)$, which they introduce differently.

The first author acknowledges support from ANR project RepRed,
ANR-09-JCJC-0102-01 and the second author acknowledges support 
from NSERC and AIM. This work began when both authors were participating
in the Combinatorial Representation Theory program at MSRI in 2008.
Hence we thank MSRI and the organizers of this program for their
hospitality and good working environment.

\section{Reflection functors}
\label{se:ReflFunc}
We begin this section by reviewing the definition of the preprojective
algebra $\Pi(Q)$ of a quiver~$Q$. Then we define and study reflection
functors for $\Pi(Q)$-modules.

\subsection{Recall on preprojective algebras}
\label{ss:PreprojAlg}
Let $K$ be a field, fixed for the whole paper. Let $Q=(I,E)$ be a
quiver with vertex set $I$ and arrow set $E$. We denote by $s$ and
$t$ the source and target maps from $E$ to $I$. We denote the path
algebra of $Q$ over $K$ by $KQ$. Then a $KQ$-module $M$ is the
data of an $I$-graded vector space $\bigoplus_{i\in I}M_i$ and of
linear maps $M_a:M_{s(a)}\to M_{t(a)}$ for each arrow $a\in E$.

To each arrow $a:i\to j$ in $E$, we associate an arrow $a^*:j\to i$
going in the opposite direction. We let $H=E\sqcup E^*$ and we
extend $*$ to $H$ by setting $(a^*)^*=a$. We extend the source and
target maps to $H$ and we set $\varepsilon(a)=1$ if $a\in E$ and
$\varepsilon(a)=-1$ if $a\in E^*$.

The quotient of the path algebra of $(I,H)$ by the ideal generated by
$$\sum_{a\in H}\varepsilon(a)aa^*$$
is called the preprojective algebra of $Q$ and is denoted by
$\Pi(Q)$. Thus a $\Pi(Q)$-module $M$ is the data of an $I$-graded
vector space $\bigoplus_{i\in I}M_i$ and of linear maps
$M_a:M_{s(a)}\to M_{t(a)}$ for each arrow $a\in H$, which
satisfy
$$\sum_{\substack{a\in H\\t(a)=i}}\varepsilon(a)M_aM_{a^*}=0$$
for each $i\in I$. In this paper, we will only consider
finite-dimensional modules.

\begin{other*}{Example}
When $Q$ has no loops nor multiple arrows, one may depict $\Pi(Q)$-modules
in a simple way, by using a symbol $i$ to represent vectors from a basis
for $M_i$ and by drawing arrows to indicate the action of the linear maps.
For instance, let $Q$ be the quiver
$$1\xleftarrow a2\xrightarrow b3$$
of type $A_3$. The $\Pi(Q)$-module $M$ represented by the diagram
$$\xymatrix@dr@=1.3em{2\ar[d]\ar[r]&3\ar[d]^{-1}\\1\ar[r]&2}$$
has dimension-vector $(\dim M_1,\dim M_2,\dim M_3)=(1,2,1)$; the map
$M_a$ sends the top basis vector in $M_2$ to the basis vector in $M_1$
and the map $M_{b^*}$ sends the basis vector in $M_3$ to the negative
of the bottom basis vector in $M_2$.
\end{other*}

The dual of a $\Pi(Q)$-module $M$ is the $\Pi(Q)$-module $M^*$
defined by taking the dual spaces and maps as follows:
$$(M^*)_i=(M_i)^*\quad\text{and}\quad(M^*)_a=(M_{a^*})^*$$
for all $i\in I$ and $a\in H$. This duality is an involutive
antiautoequivalence $*$ on the category of $\Pi(Q)$-modules.

To any $\Pi(Q)$-module $M$, one associates the tuple
$\dimvec M=(\dim M_i)$ in $\mathbb N^I$, called the
dimension-vector of $M$. Conversely, given a dimension-vector
$\nu\in\mathbb N^I$, we set $M_i=K^{\nu_i}$ and denote by
$\Lambda(\nu)$ the variety of all matrices
$$(M_a)\in\prod_{a\in H}\Hom_K(M_{s(a)},M_{t(a)})$$
that satisfy the preprojective relations
$$\sum_{\substack{a\in H\\t(a)=i}}\varepsilon(a)M_aM_{a^*}=0$$
at each vertex $i\in I$. Thus a point in $\Lambda(\nu)$ is
a representation of $\Pi(Q)$ on the vector space
$\bigoplus_{i\in I}M_i$. The group
$$G(\nu)=\prod_{i\in I}\Aut(M_i)$$
acts on $\Lambda(\nu)$ by conjugation. The orbits of this
action are in canonical bijection with the isomorphism classes
of $\Pi(Q)$-modules with dimension-vector $\nu$.

The lattice $\mathbb Z^I$ is equipped with a symmetric bilinear
form defined by
$$\langle\mu,\nu\rangle=2\sum_{i\in I}\mu_i\nu_i-\sum_{a\in H}
\nu_{s(a)}\mu_{t(a)}.$$
The following formula is due to Crawley-Boevey (\cite{Crawley-Boevey00},
Lemma~1): if $M$ and $N$ are two $\Pi(Q)$-modules, then
$$\langle\dimvec M,\dimvec N\rangle=\dim\Hom_{\Pi(Q)}(M,N)+
\dim\Hom_{\Pi(Q)}(N,M)-\dim\Ext^1_{\Pi(Q)}(M,N).$$

We denote the standard basis of $\mathbb Z^I$ by $(\alpha_i)$.
Given a vertex $i$, we denote by $S_i$ the $\Pi(Q)$-module of
dimension one concentrated at $i$ on which all arrows act as zero;
thus $\dimvec S_i=\alpha_i$. Given a $\Pi(Q)$-module $M$, the
$i$-socle of $M$, denoted $\soc_iM$, is the $S_i$-isotypic
component of the socle of $M$; likewise, the $i$-head of $M$,
denoted $\hd_iM$, is the $S_i$-isotypic component of the
head of $M$. Thus $\soc_i$ and $\hd_i$ are endofunctors on the
category of $\Pi(Q)$-modules.

\begin{other*}{Remark}
Usually, one restricts attention to representations that satisfies
a nilpotency condition (see \cite{Lusztig90b}, section~8.2). This
amounts to investigate the full subcategory of $\Pi(Q)$-mod formed
by the modules whose composition factors all belong to the set
$\{S_i\mid i\in I\}$. It can however be shown that these nilpotency
conditions are automatically satisfied in the case where $Q$ is of
type $ADE$ (see Satz~1 in \cite{Riedtmann80} or Proposition~14.2~(a)
in \cite{Lusztig91}). In this case, $\{S_i\mid i\in I\}$ is a complete
set of simple $\Pi(Q)$-modules, and thus a $\Pi(Q)$-module with
trivial $i$-socle for all $i\in I$ is itself trivial.
\end{other*}

\subsection{Definition of the reflection functors}
\label{ss:DefReflFunc}
We fix a vertex $i\in I$ at which $Q$ has no loops. Our aim is to
define a pair of adjoint endofunctors $(\Sigma_i^*,\Sigma_i)$ on
the category of $\Pi(Q)$-modules.

We break the datum of a $\Pi(Q)$-module $M$ in two parts. The first
part consists of the vector spaces $M_j$ for $j\neq i$ and of the linear
maps between them; the second part consists of the vector spaces and of
the linear maps that appear in the diagram
$$\bigoplus_{\substack{a\in H\\t(a)=i}}M_{s(a)}
\xrightarrow{(\varepsilon(a)M_a)}M_i
\xrightarrow{\ (M_{a^*})\ }\bigoplus_{\substack{a\in
H\\t(a)=i}}M_{s(a)}.$$
For brevity, we will write this diagram as
\begin{equation*}
\widetilde M_i\xrightarrow{M_{\iin(i)}}M_i
\xrightarrow{M_{\out(i)}}\widetilde M_i.\tag{$*$}
\end{equation*}
The preprojective relation at $i$ is $M_{\iin(i)}M_{\out(i)}=0$.

We will now construct a new $\Pi(Q) $-module by replacing~($*$) with
$$\widetilde M_i\xrightarrow{\overline M_{\out(i)}M_{\iin(i)}}
\ker M_{\iin(i)}\hookrightarrow\widetilde M_i,$$
where the map $\overline M_{\out(i)}:M_i\to\ker M_{\iin(i)}$ is
induced by $M_{\out(i)}$. We glue this new datum with the remaining
part of $M$. The preprojective relations will still be satisfied,
because the replacement does not change the endomorphism
$M_{\out(i)}M_{\iin(i)}$ of $\widetilde M_i$. We thus end up
with a new $\Pi(Q)$-module, which we denote by $\Sigma_iM$.

Likewise, we may replace~($*$) with
$$\widetilde M_i\twoheadrightarrow\coker M_{\out(i)}
\xrightarrow{M_{\out(i)}\overline M_{\iin(i)}}\widetilde M_i,$$
where the map $\overline M_{\iin(i)}:\coker M_{\out(i)}\to M_i$ is
induced by $M_{\iin(i)}$. We again end up with a new $\Pi(Q)$-module,
which we denote by $\Sigma_i^*M$.

One checks without difficulty that these constructions define
covariant additive functors $\Sigma_i$ and $\Sigma_i^*$, that
$\Sigma_i$ is left-exact and that $\Sigma_i^*$ is right-exact.

\begin{other*}{Example}
Let $Q$ be a quiver of type $A_3$. The diagram below presents the
action of $\Sigma_2$ and $\Sigma_2^*$ on several $\Pi(Q)$-modules.
$$
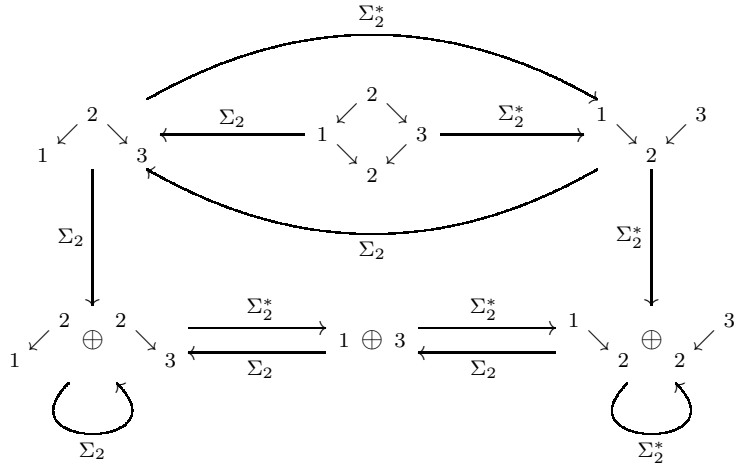
\xymatrix@=4em{&{\bsm&&2&&\\&\swarrow&&\searrow&\\1&&&&3\esm}
\ar@/^8ex/[rr]^{\Sigma_2^*}\ar[d]_{\Sigma_2}&
{\bsm&&2&&\\&\swarrow&&\searrow&\\1&&&&3\\&\searrow&&\swarrow&\\&&2&&\esm}
\ar[r]^{\Sigma_2^*}\ar[l]_{\Sigma_2}&
{\bsm1&&&&3\\&\searrow&&\swarrow&\\&&2&&\esm}
\ar@/^8ex/[ll]^{\Sigma_2}\ar[d]_{\Sigma_2^*}&\\&
{\bsm&&2\\&\swarrow&\\1&&\esm\oplus\bsm2&&\\&\searrow&\\&&3\esm}
\ar@(dl,dr)[]+/va(240)4ex/;[]+/va(300)4ex/_{\Sigma_2}
\ar@<1ex>[r]^(.6){\Sigma_2^*}&
{\bsm1\esm\oplus\bsm3\esm}
\ar@<1ex>[l]^(.4){\Sigma_2}\ar@<1ex>[r]^(.4){\Sigma_2^*}&
{\bsm1&&\\&\searrow&\\&&2\esm}\oplus{\bsm&&3\\&\swarrow&\\2&&\esm}
\ar@(dl,dr)[]+/va(240)4ex/;[]+/va(300)4ex/_{\Sigma_2^*}
\ar@<1ex>[l]^(.6){\Sigma_2}}$$
\end{other*}

\begin{other}{Remark}
\label{rk:ReflFunc}
\begin{enumerate}
\item\label{it:RFa}
We could as well distribute the signs $\varepsilon(a)$ differently in
($*$), but that would just change $\Sigma_i$ by an isomorphism. Indeed,
write $M'_{\iin(i)}$ and $M'_{\out(i)}$ for the two maps that appear in
$$ \bigoplus_{\substack{a\in H\\t(a)=i}}M_{s(a)}\xrightarrow{\ (M_a)\
}M_i\xrightarrow{(\varepsilon(a)M_{a^*})}\bigoplus_{\substack{a\in
H\\t(a)=i}}M_{s(a)}.$$
Then the module $\Sigma_iM$ defined above fits in a commutative
diagram
$$\xymatrix@C=5.5em@R=2em{
\vrule width0pt depth4pt height7pt\smash{\widetilde{(\Sigma_iM)}_i}
\ar@{=}[d]\ar[r]^{(\Sigma_iM)'_{\iin(i)}}&
(\Sigma_iM)_i\ar[d]^\simeq\ar[r]^{(\Sigma_iM)'_{\out(i)}}&
\vrule width0pt depth4pt height7pt\smash{\widetilde{(\Sigma_iM)}_i}
\ar@{=}[d]\\\widetilde M_i
\ar[r]_(.44){\overline{M'}^{}_{\!\!\out(i)}\,{M'}^{}_{\!\!\iin(i)}}&
\ker M'_{\iin(i)}\ar@{^{(}->}[r]&\widetilde M_i,}$$
where the isomorphism (middle vertical arrow) is functorial in $M$.
A similar remark can be stated about $\Sigma_i^*$. As a consequence,
we see that $\Sigma_i^*\cong*\Sigma_i*$.
\item\label{it:RFb}
The commutative diagram
$$\xymatrix@R=2em@C=5.7em{\widetilde M_i\ar@{=}[d]\ar@{->>}[r]&
\coker M_{\out(i)}\ar[d]^{\overline M_{\iin(i)}}\ar[r]^(.6)
{M_{\out(i)}\overline M_{\iin(i)}}&\widetilde M_i\ar@{=}[d]\\
\widetilde M_i\ar@{=}[d]\ar[r]^{M_{\iin(i)}}&M_i\ar[d]^{\overline
M_{\out(i)}}\ar[r]^{M_{\out(i)}}&\widetilde M_i\ar@{=}[d]\\
\widetilde M_i\ar[r]_(.43){\overline M_{\out(i)}M_{\iin(i)}}
&\ker M_{\iin(i)}\ar@{^{(}->}[r]&\widetilde M_i}$$
shows the existence of canonical morphisms
$\Sigma_i^*M\to M\to\Sigma_iM$.
\end{enumerate}
\end{other}

\begin{proposition}
\label{pr:ReflFuncAdj}
\begin{enumerate}
\item\label{it:RFAa}
The pair $(\Sigma_i^*,\Sigma_i)$ is a pair of adjoint functors.
\item\label{it:RFAb}
The adjunction morphisms $\id\to\Sigma_i\Sigma_i^*$ and
$\Sigma_i^*\Sigma_i\to\id$ can be inserted in functorial short exact
sequences
$$0\to\soc_i\to\id\to\Sigma_i\Sigma_i^*\to0\quad\text{and}\quad
0\to\Sigma_i^*\Sigma_i\to\id\to\hd_i\to0.$$
\end{enumerate}
\end{proposition}
\begin{proof}
To establish \ref{it:RFAa}, it is enough to define a pair of converse
bijections
$$\Hom_{\Pi(Q)}(M,\Sigma_iN)\cong\Hom_{\Pi(Q)}(\Sigma_i^*M,N)$$
for any modules $M$ and $N$, which are natural in $M$ and $N$.
The construction is as follows.

Consider a morphism $f:M\to\Sigma_iN$. By definition, this is a
collection of linear maps $f_j:M_j\to(\Sigma_iN)_j$, for all $j\in I$,
which intertwine the action of the arrows in $H$. Set
$$\widetilde f_i=\bigoplus_{\substack{a\in H\\t(a)=i}}f_{s(a)}:\widetilde
M_i\to\widetilde N_i.$$
In the diagram
$$\xymatrix@R=2em@C=4.7em{\widetilde M_i\ar[r]^{M_{\iin(i)}}
\ar[d]_{\widetilde f_i}&M_i\ar[r]^{M_{\out(i)}}\ar[d]^{f_i}&
\widetilde M_i\ar@{->>}[r]\ar[d]^{\widetilde f_i}&\coker
M_{\out(i)}\ar[r]^(.6){M_{\out(i)}\overline M_{\iin(i)}}
\ar@{-->}[d]^{g_i}&\widetilde M_i\ar[d]^{\widetilde f_i}\\
\widetilde N_i\ar[r]_(.43){\overline N_{\out(i)}N_{\iin(i)}}&\ker
N_{\iin(i)}\ar@{^{(}->}[r]&\widetilde N_i\ar[r]_{N_{\iin(i)}}&
N_i\ar[r]_{N_{\out(i)}}&\widetilde N_i,}$$
the two left squares commute. There is thus a unique map $g_i$ making
the third square commutative. The fourth square then also commutes,
so that if we set $g_j=f_j$ for all the vertices $j\neq i$, we get a
morphism $g:\Sigma_i^*M\to N$.

Conversely, consider a morphism $g:\Sigma_i^*M\to N$ and set
$$\widetilde g_i=\bigoplus_{\substack{a\in H\\t(a)=i}}g_{s(a)}:\widetilde
M_i\to\widetilde N_i.$$
In the diagram
$$\xymatrix@R=2em@C=4.7em{\widetilde M_i\ar[r]^{M_{\iin(i)}}
\ar[d]_{\widetilde g_i}&M_i\ar[r]^{M_{\out(i)}}\ar@{-->}[d]^{f_i}&
\widetilde M_i\ar@{->>}[r]\ar[d]^{\widetilde g_i}&\coker
M_{\out(i)}\ar[r]^(.6){M_{\out(i)}\overline M_{\iin(i)}}
\ar[d]^{g_i}&\widetilde M_i\ar[d]^{\widetilde g_i}\\
\widetilde N_i\ar[r]_(.43){\overline N_{\out(i)}N_{\iin(i)}}&\ker
N_{\iin(i)}\ar@{^{(}->}[r]&\widetilde N_i\ar[r]_{N_{\iin(i)}}&
N_i\ar[r]_{N_{\out(i)}}&\widetilde N_i,}$$
the two right squares commute. There is thus a unique map $f_i$ making
the second square commutative. The first square then also commutes,
so that if we set $f_j=g_j$ for all the vertices $j\neq i$, we get a
morphism $f:M\to\Sigma_iN$.

To establish \ref{it:RFAb}, one checks that $\Sigma_i^*\Sigma_iM$
is the $\Pi(Q)$-module obtained by replacing in $M$ the part summed up
by ($*$) with
$$\widetilde M_i\xrightarrow{M_{\iin(i)}}\im M_{\iin(i)}
\xrightarrow{M_{\out(i)}}\widetilde M_i$$
and that $\Sigma_i\Sigma_i^*M$ is the $\Pi(Q)$-module obtained by
replacing in $M$ the part summed up by ($*$) with
$$\widetilde M_i\xrightarrow{M_{\iin(i)}}\im M_{\out(i)}
\xrightarrow{M_{\out(i)}}\widetilde M_i.$$
It remains to observe that as vector spaces,
$\hd_iM\cong\coker M_{\iin(i)}$ and $\soc_iM\cong\ker M_{\out(i)}$.
\end{proof}

Proposition~\ref{pr:ReflFuncAdj}~\ref{it:RFAb} implies readily that
$\Sigma_i$ and $\Sigma_i^*$ define inverse equivalence of categories
$$\left\{\begin{gathered}
\text{$\Pi(Q)$-modules}\\\text{with trivial
$i$-head}\end{gathered}\right\}
\xymatrix{\ar@<1ex>[r]^{\Sigma_i}&\ar@<1ex>[l]^{\Sigma_i^*}}
\left\{\begin{gathered}\text{
$\Pi(Q)$-modules}\\\text{with trivial
$i$-socle}\end{gathered}\right\}$$
and that there are natural isomorphisms
$$\Sigma_i\Sigma_i^*\Sigma_i\cong\Sigma_i\quad\text{and}\quad
\Sigma_i^*\Sigma_i\Sigma_i^*\cong\Sigma_i^*.$$

Let $s_i$ be the reflection $\mu\mapsto\mu-\langle\alpha_i,\mu\rangle
\alpha_i$ on the lattice $\mathbb Z^I$. Routine arguments show that
$$\hd_iM=0\ \Longrightarrow\ \dimvec\Sigma_iM=s_i(\dimvec M).$$
We leave as an exercise for the reader to prove that
$$\hd_iM=0\ \Longleftrightarrow\ \dimvec\Sigma_iM=s_i(\dimvec M)\
\Longleftrightarrow\ M\cong\Sigma_i^*\Sigma_iM$$
and
$$\soc_iM=0\ \Longleftrightarrow\ \dimvec\Sigma_i^*M=s_i(\dimvec M)\
\Longleftrightarrow\ M\cong\Sigma_i\Sigma_i^*M;$$
we will however not use this more complete result.

\subsection{Further properties}
\label{ss:FurtProp}
We begin this section with an easy lemma, which is used in the
proof of Proposition~\ref{pr:NGammaFoncRefl}.

\begin{lemma}
\label{le:CaseTrivISoc}
If the $\Pi(Q)$-module $M$ has trivial $i$-socle, then the canonical
morphism $M\to\Sigma_iM$ is a monomorphism which induces isomorphisms
$\Sigma_iM\cong\Sigma_i^2M$ and $\soc M\cong\soc(\Sigma_iM)$.
\end{lemma}
\begin{proof}
The $i$-socle of $M$ is the kernel of the map $M_{\out(i)}$. Assume
that it is trivial. Then $\overline M_{\out(i)}$ is injective. The
canonical morphism $f:M\to\Sigma_iM$ is thus a monomorphism, and the
equality $\ker M_{\iin(i)}=\ker\left(\overline
M_{\out(i)}M_{\iin(i)}\right)$ holds, which implies
$\Sigma_iM\cong\Sigma_i^2M$.

Now consider a monomorphism $g:T\to\Sigma_iM$ from a completely
reducible $\Pi(Q)$-module $T$ into $\Sigma_iM$. Since $\Sigma_iM$
has trivial $i$-socle, $T$ cannot contain a simple submodule
isomorphic to $S_i$, hence cannot map non-trivially to $S_i$.
However by construction, the cokernel of $f$ is a direct sum of
copies of $S_i$. We conclude that $g$ factorizes through $f$.
This proves that $f$ maps $\soc M$ onto $\soc\Sigma_iM$.
\end{proof}

In~\cite{Crawley-Boevey-Holland98}, Crawley-Boevey and Holland
define reflection functors on the category of modules over the
\textit{deformed\/} preprojective algebra. Their functors define
an action of the Weyl group. Our functors do not
enjoy the property $\Sigma_i^2=\id$, but it is nevertheless
reassuring to see that they satisfy the braid relations.

\begin{proposition}
\label{pr:BraidRel}
Let $i$ and $j$ be two vertices that are linked by a single edge
in $Q$. Then the functors $\Sigma_i\Sigma_j\Sigma_i$ and
$\Sigma_j\Sigma_i\Sigma_j$ are isomorphic.
\end{proposition}
\begin{proof}
We only sketch the proof, since we will not use this property later on.
Let $(c,c^*)$ be the pair of arrows in $H$ between the vertices $i$ and
$j$. Let $M$ be a $\Pi(Q)$-module. Abbreviate the part
\vspace{-4ex}
$$\xymatrix@C=5em{
\raisebox{-7ex}{$\displaystyle\bigoplus_{\substack{a\in
H\\t(a)=i\\a\neq c^*}}M_{s(a)}$}
\ar@<.5ex>[r]^(.6){(\varepsilon(a)M_a)}&
M_i\ar@<.5ex>[l]^(.4){(M_{a^*})}\ar@<.5ex>[r]^{M_c}&
M_j\ar@<.5ex>[l]^{M_{c^*}}\ar@<.5ex>[r]^(.4){(M_{b^*})}&
\raisebox{-7ex}{$\displaystyle\bigoplus_{\substack{b\in
H\\t(b)=j\\b\neq c}}M_{s(b)}$}
\ar@<.5ex>[l]^(.6){(\varepsilon(b)M_b)}}$$
of the datum of $M$ by the notation
$$\xymatrix{R\ar@<.5ex>[r]^k&V\ar@<.5ex>[l]^l\ar@<.5ex>[r]^f&
W\ar@<.5ex>[l]^g\ar@<.5ex>[r]^m&S.\ar@<.5ex>[l]^n}$$
The preprojective relations at $i$ and $j$ then read
$$kl+\varepsilon(c^*)gf=0\quad\text{ and }\quad
nm+\varepsilon(c)fg=0.$$
Explicit computations show that applying $\Sigma_i\Sigma_j\Sigma_i$
or $\Sigma_j\Sigma_i\Sigma_j$ to $M$ both amount to replacing this
part of the datum of $M$ with
$$\xymatrix@C=5.5em{R\ar@<.5ex>[r]^{\left(\bsm lk\\mfk\esm\right)}&
V'\ar@<.5ex>[l]^{(\bsm1&\ 0\esm)}
\ar@<.5ex>[r]^{\left(\bsm0\ &\ 1\\\varepsilon(c)lk\ &\ 0\esm\right)}&
W'\ar@<.5ex>[r]^{(\bsm1&\ 0\esm)}
\ar@<.5ex>[l]^{\left(\bsm0&\ 1\\\varepsilon(c^*)mn&\ 0\esm\right)}&
S,\ar@<.5ex>[l]^{\left(\bsm mn\\lgn\esm\right)}}$$
where $V'\subseteq R\oplus S$ and $W'\subseteq S\oplus R$ are the
kernels of the linear maps
$$R\oplus S\xrightarrow{(\bsm fk&\;n\esm)}W\quad\text{ and }\quad
S\oplus R\xrightarrow{(\bsm gn&\;k\esm)}V.$$
\end{proof}

\section{The modules $N(\gamma)$}
\label{se:ModNGamma}
From now on, we fix a Dynkin diagram $\Gamma$ of type $ADE$. Orienting
the edges of this diagram yields a quiver $Q=(I,E)$. The datum of the
diagram $\Gamma$ is equivalent to that of the isomorphism class of a
simply-laced semisimple complex Lie algebra $\mathfrak g$. The set
$I$ of vertices of $\Gamma$ indexes the simple roots $\alpha_i$ of
$\mathfrak g$. The root lattice $Q$ is then the lattice $\mathbb Z^I$
with basis $(\alpha_i)_{i\in I}$; we write $Q_+$ for the image of
$\mathbb N^I$ in this identification. (We are here using the same
notation $Q$ for the quiver and for the root lattice. These
notations are traditional and we are convinced that our choice will
not be a source of confusion for the reader.) Since we are in the simply
laced case, there is no need to distinguish between roots and coroots.
The weight lattice is thus the dual $P$ of this lattice $Q$; it is
endowed with the basis $(\omega_i)_{i\in I}$ of fundamental weights.
By duality, the reflections $s_i$ defined at the end of
section~\ref{ss:DefReflFunc} act on $P$ and generate the Weyl group
$W\subseteq\Aut(P)$. Then $(W,(s_i)_{i\in I})$ is a finite Coxeter
system; we denote its length function by $\ell$ and its longest
element by $w_0$. Finally, the bilinear form $\langle\;,\;\rangle$
from section \ref{ss:PreprojAlg} is the standardly normalised
$W$-invariant scalar product on $Q$; it embeds $Q$ as a sublattice of $P$.

Our aim in this section is to define a family of $\Pi(Q)$-modules
$N(\gamma)$ indexed by weights $\gamma\in P$. These modules can be
characterized up to isomorphism in the following fashion.

\begin{theorem}
\label{th:CharNGamma}
\begin{enumerate}
\item\label{it:CNGa}
If $\gamma$ is antidominant, then $N(\gamma)=0$.
\item\label{it:CNGb}
Let $i\in I$ and let $\gamma$ be a $W$-conjugate of $-\omega_i$,
with $\gamma\neq-\omega_i$. Then $N(\gamma)$ satisfies
$$\dimvec N(\gamma)=\gamma+\omega_i\quad\text{ and }\quad
\soc N(\gamma)\cong S_i.$$
Moreover these conditions characterize $N(\gamma)$ up to
isomorphism.
\item\label{it:CNGc}
If $\gamma$ and $\delta$ belong to the same Weyl chamber, then
$N(\gamma+\delta)\cong N(\gamma)\oplus N(\delta)$.
\end{enumerate}
\end{theorem}

This theorem will be proved in section~\ref{ss:ModNGamma}.

\subsection{Nakajima's quiver}
\label{ss:NakaQuiv}
The easiest way to define the modules $N(\gamma)$ uses Nakajima's
trick (see \cite{Nakajima94}) of expanding $Q$ by adding an extra
vertex $i'$ and an extra arrow $d_i:i\to i'$ for each vertex $i\in I$.
Playing the same game as in section~\ref{ss:PreprojAlg}, we double all
arrows in this extended quiver, obtaining thereby a set of arrows
$\widehat H=H\sqcup\{d_i,d_i^*\mid i\in I\}$ with an involution $*$.
Here is an example depicting the situation for the type $A_3$.
$$\xymatrix@C=5em@R=2em{1\ar@<-.5ex>[d]_{d_1}\ar@<-.5ex>[r]_{a^*}&2
\ar@<-.5ex>[d]_{d_2}\ar@<-.5ex>[l]_a\ar@<.5ex>[r]^b&3
\ar@<-.5ex>[d]_{d_3}\ar@<.5ex>[l]^{b^*}\\1'\ar@<-.5ex>[u]_{d_1^*}
&2'\ar@<-.5ex>[u]_{d_2^*}&3'.\ar@<-.5ex>[u]_{d_3^*}}$$

We then construct the preprojective algebra of the extended quiver,
which we denote by $\Pi(\widehat Q)$. Hence a $\Pi(\widehat Q)$-module
consists of two $I$-graded vector spaces $\bigoplus_{i\in I}M_i$ and
$\bigoplus_{i\in I}M_{i'}$ along with linear maps $M_a:M_{s(a)}\to
M_{t(a)}$, for $a\in H$, and $M_{d_i}:M_i\to M_{i'}$ and
$ M_{d^*_i}:M_{i'}\to M_i$, for $i\in I$, which satisfy the
preprojective relations
$$\sum_{\substack{a\in H\\t(a)=i}}\varepsilon(a)M_aM_{a^*}=
M_{d_i^*}M_{d_i}\quad\text{ and }\quad M_{d_i}M_{d_i^*}=0$$
at each vertex $i\in I$. The dimension-vector of $M$ is defined as
the pair $\dimvec M=(\nu,\lambda)$ in $\mathbb N^I\times\mathbb N^{I'}$,
where $\nu_i=\dim M_i$ and $\lambda_{i'}=\dim M_{i'}$.

Given a dimension-vector $(\nu,\lambda)\in\mathbb N^I\times\mathbb
N^{I'}$ for $\Pi(\widehat Q)$, it is customary to view
$\lambda=(\lambda_{i'})$ as the antidominant weight $-\sum_{i\in I}
\lambda_{i'}\omega_i$ in $P$. The main justification for this
identification is that the pairing between $P$ and $Q$ coincides
with the restriction to $\mathbb Z^{I'}$ and $\mathbb Z^I$ of the
symmetric bilinear form of the extended quiver. This fact will be
exploited below when we write expressions like
$$\langle(\nu,\lambda),(\alpha_i,0)\rangle=\langle\nu+\lambda,
\alpha_i\rangle;$$
here the left hand side is the symmetric bilinear form from
section~\ref{ss:PreprojAlg} and the right hand side is the
duality pairing between $P$ and $Q$.

\subsection{Stable $\Pi(\widehat Q)$-modules}
\label{ss:StabMod}
We say that a $\Pi(\widehat Q)$-module $M$ is stable if for each
$i\in I$, the linear map
$$M_i\xrightarrow{\left(\bsm(M_a)\\M_{d_i}\esm\right)}\left(
\bigoplus_{\substack{a\in H\\s(a)=i}}M_{t(a)}\right)\oplus M_{i'}$$
is injective. In other words, we ask that $\soc_iM=0$ at each
non-primed vertex.

This notion of stability was introduced by Nakajima. Indeed using
the remark in section \ref{ss:PreprojAlg}, one can easily see that
this definition is equivalent to Definition~3.9 in \cite{Nakajima98}.
We now study the stable $\Pi(\widehat Q)$-modules following ideas of
Nakajima and Saito.

\begin{lemma}
\label{le:HeadBound}
Let $M$ be a stable $\Pi(\widehat{Q})$-module with dimension-vector $(\nu,\lambda)$.
Let $c$ be the dimension of the $i$-head of $M$. Then
$$\dim \Hom_{\Pi(\widehat Q)}(M,S_i)=c\quad\text{and}\quad
\dim\Ext^1_{\Pi(\widehat Q)}(M,S_i)=c-\langle\nu+\lambda,\alpha_i\rangle.$$
In particular, $c\geq\max(0,\langle\nu+\lambda,\alpha_i\rangle)$.
\end{lemma}
\begin{proof}
The first equality is equivalent to the assertion $\dim\hd_iM=c$.
On the other hand, the stability condition means that
$\Hom_{\Pi(\widehat Q)}(S_i,M)=0$. Now the dimension-vector of $S_i$
is $(\alpha_i,0)$, hence the result follows from Crawley-Boevey's
formula (see section~\ref{ss:PreprojAlg}).
\end{proof}

Let $M$ be a stable $\Pi(\widehat{Q})$-module with dimension-vector
$(\nu,\lambda)$. By Lemma~\ref{le:HeadBound}, the $i$-head of $M$
has dimension at least $\max(0,\langle\nu+\lambda,\alpha_i\rangle)$.
If there is equality here, then we say that $M$ has a small $i$-head.

\begin{lemma}
\label{le:AnalRep}
\begin{enumerate}
\item\label{it:ARa}
Let $M$ be a stable $\Pi(\widehat Q)$-module with dimension-vector $(\nu,\lambda)$
whose $i$-head has dimension $c$. Then for any integer $k$ such that
$-c\leq k\leq c-\langle\nu+\lambda,\alpha_i\rangle$, there exists a
stable $\Pi(\widehat Q)$-module with dimension-vector $(\nu+k\alpha_i,\lambda)$.
\item\label{it:ARb}
Assume that $\langle\nu+\lambda,\alpha_i\rangle\leq0$. Then the functors
$\Sigma_i$ and $\Sigma_i^* $ restrict to equivalences of categories
\begin{equation*}
\left\{\begin{gathered}
\text{stable $\Pi(\widehat Q)$-modules with}\\
\text{dimension-vector $(\nu,\lambda)$}\\
\text{with small $i$-head}
\end{gathered}\right\}
\xymatrix{\ar@<1ex>[r]^{\Sigma_i}&\ar@<1ex>[l]^{\Sigma_i^*}}
\left\{\begin{gathered}
\text{stable $\Pi(\widehat Q)$-modules with}\\
\text{dimension-vector $(\nu-\langle\nu+\lambda,\alpha_i\rangle,\lambda)$}\\
\text{with small $i$-head}
\end{gathered}\right\}.
\end{equation*}
\item\label{it:ARc}
If there exists a stable $\Pi(\widehat Q)$-module with dimension-vector
$(\nu,\lambda)$, then
$$\nu+\lambda\in\bigcap_{w\in W}w(\lambda+Q_+).$$
\end{enumerate}
\end{lemma}

\begin{proof}
Items~\ref{it:ARa} and~\ref{it:ARb} are Lemma~4.2.1 and a special
case of Lemma~4.2.2 in~\cite{Saito02}. Item~\ref{it:ARc} is a consequence
of Corollary~10.8 in~\cite{Nakajima98}. We nevertheless recall the proof
of all these facts for the convenience of the reader.

Let $M$ be a stable $\Pi(\widehat Q)$-module with dimension-vector 
$(\nu,\lambda)$.
We analyse the situation locally around the vertex $i$. For brevity,
let us introduce the notation $\widetilde M_i$, $M_{\iin(i)}$ and
$M_{\out(i)}$ for the objects in the diagram
$$\left(\bigoplus_{\substack{a\in H\\t(a)=i}}M_{s(a)}\right)\oplus M_{i'}
\xrightarrow{\left(\bsm(\varepsilon(a)M_a)&-M_{d_i^*}\esm\right)}M_i
\xrightarrow{\quad\left(\bsm(M_{a^*})\\M_{d_i}\esm\right)\quad}
\left(\bigoplus_{\substack{a\in H\\t(a)=i}}M_{s(a)}\right)\oplus M_{i'},$$
as we did in section~\ref{ss:DefReflFunc}. The stability condition means
that $M_{\out(i)}$ is injective.

Consider the endomorphism $u=M_{\out(i)}M_{\iin(i)}$ of
$\widetilde M_i$. Note that
$$\im u\subseteq\im M_{\out(i)}\subseteq\ker u,$$
that
$$\dim \im M_{\out(i)}-\dim\im u=\dim\coker M_{\iin(i)}=c,$$
and that
\begin{align*}
\dim\ker u-\dim\im M_{\out(i)}
&=\dim\widetilde M_i-\dim\im u-\dim\im M_{\out(i)}\\
&=c+\dim\widetilde M_i-2\dim M_i\\
&=c-\langle\nu+\lambda,\alpha_i\rangle.
\end{align*}

Thus for any $k$ satisfying $-c\leq k\leq c-\langle\nu+\lambda,
\alpha_i\rangle$, we can find a subspace $V\subseteq\widetilde M_i$
of dimension $k+\dim M_i$ satisfying $\im u\subseteq V\subseteq\ker u$.
We can then construct a new $\Pi(\widehat Q)$-module by replacing
$$\widetilde M_i\xrightarrow{M_{\iin(i)}}M_i
\xrightarrow{M_{\out(i)}}\widetilde M_i\qquad\text{with}\qquad
\widetilde M_i\xrightarrow{\;u\;}V\hookrightarrow\widetilde M_i$$
in $M$. We thus get a stable $\Pi(\widehat Q)$-module of
dimension-vector $(\nu+k\alpha_i,\lambda)$, which
establishes~\ref{it:ARa}.

Keeping the same notation, we now observe that $M$ has a small $i$-head
if and only if either $c$ or $c-\langle\nu+\lambda,\alpha_i\rangle$
equals zero. In other words, $M$ has a small $i$-head if and only if
$\im M_{\out(i)}$ is either $\im u$ or $\ker u$. Hence keeping
$\widetilde M_i$ and $u$ the same and flipping $M_i$ between $\im u$
and $\ker u$ yields equivalences. This flipping is indeed what the
functors $\Sigma_i$ and $\Sigma_i^*$ do: under the assumption
$\langle\nu+\lambda,\alpha_i\rangle\leq0$, if $M$ has dimension-vector
$(\nu,\lambda)$, then we deal with the first possibility $\im
M_{\out(i)}=\im u$, and $\Sigma_i M$ is the same module except with
$M_i$ replaced by $\ker M_{\iin(i)}=\ker u$. And if $M$ has dimension-vector
$(\nu-\langle\nu+\lambda,\alpha_i\rangle\alpha_i,\lambda)$, then
we deal with the second possibility $\im M_{\out(i)}=\ker u$, and
then $\Sigma_i^* M$ is the same module, except with $M_i$ replaced
by $\coker M_{\out(i)}=\widetilde M_i/\ker u\cong\im u$. We have thus
established~\ref{it:ARb}.

Looking now at~\ref{it:ARc}, we fix $\lambda$ and allow $\nu$ to vary.
Assertion~\ref{it:ARa} implies that the set
$$\bigl\{\nu+\lambda\in P\bigm|\text{there exists a stable
$\Pi(\widehat Q)$-module with dimension-vector $(\nu,\lambda)$}\bigr\}$$
is $W$-invariant. Since it is a subset of $\lambda+Q_+$, it is
contained in $\bigcap_{w\in W}w(\lambda+Q_+)$. This
establishes~\ref{it:ARc}.
\end{proof}

\subsection{The modules $\widehat N(\gamma)$}
\label{ss:ModHatNGamma}
Lemma~\ref{le:AnalRep} allows us to quickly prove the following
theorem, which can be seen as a corollary to the work of
Nakajima~\cite{Nakajima98} and Lusztig~\cite{Lusztig98}.
\begin{theorem}
\label{th:ModHatNGamma}
Let $\lambda$ be an antidominant weight and let $w\in W$. Then there
exists a unique (up to isomorphism) stable $\Pi(\widehat Q)$-module
with dimension-vector $(w\lambda-\lambda,\lambda)$. This module has a small
$i$-head, for all $i\in I$.
\end{theorem}
\begin{proof}
Let $\lambda$ and $w$ be as in the statement of the theorem.

We first prove the second assertion. Let $M$ be a stable
$\Pi(\widehat Q)$-module $M$ with dimension-vector
$(w\lambda-\lambda,\lambda)$ and let $i\in I$. We set
$c=\dim\hd_iM$ and distinguish two cases.

If $\ell(s_iw)>\ell(w)$, then we take $k=-c$ in Lemma~\ref{le:AnalRep}
\ref{it:ARa} and get the existence of a stable $\Pi(\widehat Q)$-module
with dimension-vector $(w\lambda-\lambda-c\alpha_i,\lambda)$.
Lemma~\ref{le:AnalRep}~\ref{it:ARc} then says that
$w\lambda-c\alpha_i\in w(\lambda+Q_+)$, and therefore
$-cw^{-1}(\alpha_i)\in Q_+$. Since $w^{-1}(\alpha_i)$ is a positive
root here, we conclude that $c=0$. Moreover $\langle w\lambda,\alpha_i
\rangle\leq0$ by the antidominance of $\lambda$.

If $\ell(s_iw)<\ell(w)$, then we take $k=c-\langle w\lambda,\alpha_i\rangle$
in Lemma~\ref{le:AnalRep}~\ref{it:ARa} and get the existence of a
stable $\Pi(\widehat Q)$-module with dimension-vector
$(w\lambda-\lambda+k\alpha_i,\lambda)$. Lemma~\ref{le:AnalRep}~\ref{it:ARc}
then says that $w\lambda+k\alpha_i\in w(\lambda+Q_+)$. Since
$w^{-1}(\alpha_i)$ is a negative root here, we obtain $k\leq0$. But
also $k\geq0$, by Lemma~\ref{le:HeadBound}, and therefore
$c=\langle w\lambda,\alpha_i\rangle$.

So in both cases, $c=\max(0,\langle w\lambda,\alpha_i\rangle)$, as claimed.

We now turn to the uniqueness. The result is obvious for $w=1$.
If $w\neq1$, we pick $i\in I$ such that $\ell(s_iw)<\ell(w)$. Then
Lemma~\ref{le:AnalRep}~\ref{it:ARb} ensures that
\begin{equation*}
\left\{\begin{gathered}
\text{stable $\Pi(\widehat Q)$-modules with}\\
\text{dimension-vector $(s_iw\lambda-\lambda,\lambda)$}
\end{gathered}\right\}
\xymatrix{\ar@<1ex>[r]^{\Sigma_i}&\ar@<1ex>[l]^{\Sigma_i^*}}
\left\{\begin{gathered}
\text{stable $\Pi(\widehat Q)$-modules with}\\
\text{dimension-vector $(w\lambda-\lambda,\lambda)$}
\end{gathered}\right\}
\end{equation*}
is an equivalence of categories, for the condition about the smallness
of the $i$-head is automatically fulfilled. Thus the uniqueness result
for $w$ is equivalent to that for $s_iw$. At this point, an induction
concludes the proof.
\end{proof}

Since the datum of $w\lambda$ allows to recover $\lambda$, we may
denote by $\widehat N(w\lambda)$ the $\Pi(\widehat Q)$-module whose
existence and uniqueness is asserted by Theorem~\ref{th:ModHatNGamma}.
The following easy proposition study the behaviour of these modules
under the reflection functors.

\begin{proposition}
\label{pr:NGammaFoncRefl}
Let $\gamma$ be a weight and let $i\in I$. If
$\langle\gamma,\alpha_i\rangle\leq0$, then
$$\Sigma_i\widehat N(\gamma)\cong\widehat N(s_i\gamma),\quad
\Sigma_i\widehat N(s_i\gamma)\cong\widehat N(s_i\gamma)\quad
\text{ and }\quad
\Sigma_i^*\widehat N(s_i\gamma)\cong\widehat N(\gamma).$$
\end{proposition}
\begin{proof}
The first relation is a direct consequence of the proof of Theorem
\ref{th:ModHatNGamma}; to see it, it suffices to write
$\gamma=w\lambda$ with $\lambda$ antidominant and $w\in W$ such
that $\ell(s_iw)>\ell(w)$.

Lemma~\ref{le:CaseTrivISoc} then gives the second relation:
$$\Sigma_i\widehat N(s_i\gamma)\cong\Sigma_i^2\widehat N(\gamma)
\cong\Sigma_i\widehat N(\gamma)\cong\widehat N(s_i\gamma).$$

Finally the third relation comes from Proposition~\ref{pr:ReflFuncAdj}
and from the fact that $\widehat N(\gamma)$ has trivial $i$-head:
$$\Sigma_i^*\widehat N(s_i\gamma)\cong\Sigma_i^*\Sigma_i\widehat
N(\gamma)\cong\widehat N(\gamma).$$
\end{proof}

To conclude this section, we record the following consequence of
Theorem~\ref{th:ModHatNGamma} and Lemma~\ref{le:HeadBound}: for
any weight $\gamma$ and any $i\in I$,
\begin{equation}
\left\{\begin{aligned}
&\dim\Hom_{\Pi(\widehat Q)}\bigl(\widehat N(\gamma),S_i\bigr)
=\max(0,\langle\gamma,\alpha_i)\rangle,\\[4pt]
&\dim\Ext^1_{\Pi(\widehat Q)}\bigl(\widehat N(\gamma),S_i\bigr)
=\max(0,-\langle\gamma,\alpha_i\rangle).
\end{aligned}\right.
\label{eq:DimHomNGamma}
\end{equation}

\subsection{The modules $N(\gamma)$}
\label{ss:ModNGamma}
Let $\lambda$ be an antidominant weight and let $w\in W$.
Looking at the proofs of Lemma~\ref{le:AnalRep}~\ref{it:ARb} and
Theorem~\ref{th:ModHatNGamma}, one easily proves that if
$\ell(s_iw)>\ell(w)$, then there is a canonical monomorphism
$\widehat N(w\lambda)\hookrightarrow\widehat N(s_iw\lambda)$.
(Alternatively, one can use Lemma~\ref{le:CaseTrivISoc} and
Proposition~\ref{pr:NGammaFoncRefl}.) In particular,
$\widehat N(\lambda)$ is contained in each $\widehat N(w\lambda)$.
Examining the dimension-vectors, one sees that $\widehat N(\lambda)$
is the sum at the prime vertices of the vector spaces of
$\widehat N(w\lambda)$; the inclusion $\widehat N(\lambda)
\hookrightarrow\widehat N(w\lambda)$ therefore means that the arrows
$d_i^*$ act as the zero map on the module $\widehat N(w\lambda)$.
We set $N(w\lambda)=\widehat N(w\lambda)/\widehat N(\lambda)$.

Now a $\Pi(Q)$-module $T$ can be viewed as a
$\Pi(\widehat Q)$-module by setting $T_{i'}=0$ for all $i\in I$.
Conversely, the $\Pi(\widehat Q)$-module $N(w\lambda)$ can be
viewed as a $\Pi(Q)$-module, because it is supported only on
non-primed vertices. In this context, there are natural isomorphisms
$$\Hom_{\Pi(\widehat Q)}\bigl(\widehat N(w\lambda),T\bigr)\cong
\Hom_{\Pi(\widehat Q)}\bigl(N(w\lambda),T\bigr)\cong
\Hom_{\Pi(Q)}\bigl(N(w\lambda),T\bigr).$$

We are now in a position to prove Theorem~\ref{th:CharNGamma}.
\trivlist
\item[\hskip\labelsep{\itshape Proof of Theorem~\ref{th:CharNGamma}.}]
\upshape
Assertion~\ref{it:CNGa} is trivial.

Consider now the assertion~\ref{it:CNGb} and write $\gamma=-w\omega_i$,
with $i\in I$ and $w\in W$. Set $M=\widehat N(\gamma)$; this is a stable
$\Pi(\widehat Q)$-module with dimension-vector $(\gamma+\omega_i,-\omega_i)$.
Then the vector space $\bigoplus_{j\in I}M_{j'}$ coincides with the module
$\widehat N(-\omega_i)$, so has dimension one and is concentrated at
vertex~$i'$. On the other hand, the $\Pi(Q)$-module $N(\gamma)$ identifies
with the vector space $\bigoplus_{i\in I}M_i$. In particular, the
dimension-vector of $N(\gamma)$ is $\gamma+\omega_i$. For each $j\in I$,
the $j$-socle of $M$ is the intersection of the $j$-socle of $N(\gamma)$
with $\ker M_{d_j}$. The stability of $M$ therefore implies that the
socle of $N(\gamma)$ is concentrated at vertex $i$ and at most
one-dimensional. Since $N(\gamma)\neq0$, this socle cannot be zero, and
we conclude that $\soc N(\gamma)\cong S_i$. The module $N(\gamma)$ thus
enjoys the properties stated in assertion~\ref{it:CNGb}. Conversely,
let $N$ be a $\Pi(Q)$-module that enjoys these properties. One can
then construct a stable $\Pi(\widehat Q)$-module $M$ with dimension-vector
$(\gamma+\omega_i,-\omega_i)$ as follows: $M_j=N_j$ for all $j\in I$,
$M_h=N_h$ for all $h\in H$, $M_{i'}=\soc_iN$, $M_{d_i}:N_i\to\soc_iN$ is
any linear map which extends the identity of $\soc_iN$, and the
remaining spaces and arrows are zero. By the uniqueness in
Theorem~\ref{th:ModHatNGamma}, $M$ is isomorphic to $\widehat N(\gamma)$,
and we conclude that $N$ is isomorphic to $N(\gamma)$.
Assertion~\ref{it:CNGb} is proved.

Finally, let $\lambda$ and $\mu$ be two antidominant weights and let
$w\in W$. Then we have an isomorphism $\widehat N(w\lambda)\oplus\widehat
N(w\mu)\cong\widehat N(w(\lambda+\mu))$, again by the uniqueness in
Theorem~\ref{th:ModHatNGamma}. From there, it follows easily that
$N(w\lambda)\oplus N(w\mu)\cong N(w(\lambda+\mu))$, which is the
content of~\ref{it:CNGc}.
\nobreak\noindent$\square$\endtrivlist

The modules $N(\gamma)$ have already been studied at great length by
Gei\ss, Leclerc and Schr\"oer, who have shown their relevance to the
theory of cluster algebras. More precisely, the reflection functors
$\Sigma_i$ can be seen as partial inverses to the functors
$\mathcal E_i$ defined in section 5 of~\cite{Geiss-Leclerc-Schroer08}.
The key to this interpretation is Proposition~\ref{pr:NGammaFoncRefl}
and the following result.

\begin{proposition}
\label{pr:ProjCover}
The projective cover of $S_i$ in $\Pi(Q)$-mod is $N(\omega_i)$.
\end{proposition}
\begin{proof}
Consider a simple $\Pi(Q)$-module $S_j$, with $j\in I$. In the
exact sequence
$$\Hom_{\Pi(\widehat Q)}\bigl(\widehat N(w_0\omega_i),S_j\bigr)\to
\Ext^1_{\Pi(\widehat Q)}\bigl(N(\omega_i),S_j\bigr)\to
\Ext^1_{\Pi(\widehat Q)}\bigl(\widehat N(\omega_i),S_j\bigr),$$
the two extreme terms are zero, thanks to the equations
(\ref{eq:DimHomNGamma}). Hence the middle term is also zero. Observing
that a non-split extension of $N(\omega_i)$ by $S_j$ in the category
$\Pi(Q)$-mod would also be non-split in $\Pi(\widehat Q)$-mod, we see
that $\Ext^1_{\Pi(Q)}(N(\omega_i),S_j)=0$. We conclude that
$\Ext^1_{\Pi(Q)}(N(\omega_i),T)=0$ for any object $T$ in the
category $\Pi(Q)$-mod, because such a $T$ has always a finite
filtration with subquotients $S_j$.

Thus $N(\omega_i)$ is a projective $\Pi(Q)$-module, that is, is the
projective cover of its head. On the other hand, the equality
$$\dim\Hom_{\Pi(Q)}(N(\omega_i),S_j)=\dim\Hom_{\Pi(\widehat
Q)}(\widehat N(\omega_i),S_j)=\langle\alpha_j,\omega_i\rangle$$
shows that the head of $N(\omega_i)$ is $S_i$.
\end{proof}

We conclude this section by noting that in \cite{Kamnitzer-Sadanand10},
Sadanand and the second author present a complete and explicit
description of the modules $N(\gamma)$ in type A.

\section{Pseudo-Weyl polytopes}
\label{se:PWPoly}
For any weight $\gamma$ and any $\Pi(Q)$-module $T$, we set
$D_\gamma(T)=\dim\Hom_{\Pi(Q)}(N(\gamma),T)$. Thus
$D_\gamma(T)$ is always non-negative and $D_\gamma(T)=0$
if $\gamma$ is antidominant.

In this section, we study these functions $D_\gamma$. We first note
the following consequence of Theorem~\ref{th:CharNGamma}~\ref{it:CNGc}:
$D_{\gamma+\delta}=D_\gamma+D_\delta$ whenever $\gamma$ and
$\delta$ belong to the same Weyl chamber. This observation prompts
us to pay a particular attention to the so-called chamber weights,
that is, weights of the form $\gamma=w\omega_i$, with $w\in W$ and
$i\in I$.

\subsection{Edge and other relations}
\label{ss:EdgeOthRel}
\begin{proposition}
\label{pr:BZDataPrep}
Let $\gamma$ be a weight, let $i\in I$ and let $T$ be a $\Pi(Q)$-module.
If $\langle\gamma,\alpha_i\rangle\leq0$, then
\begin{align}
D_{\gamma}(T)&=D_{s_i\gamma}(\Sigma_iT),
\label{eq:BZDataRefl}\\
D_{s_i\gamma}(T)&=D_{\gamma}(\Sigma_i^*T)-
\langle\gamma,\dimvec\soc_iT\rangle,
\label{eq:InductRelPrep}\\
D_{\gamma}(T)&=D_{\gamma}(\Sigma_i^*\Sigma_iT).
\label{eq:BZDataRadi}
\end{align}
\end{proposition}
\begin{proof}
We can regard $T$ as a module over either $\Pi(Q)$ or $\Pi(\widehat Q)$;
this does not affect the formation of $\Sigma_iT$. Having noticed this,
equation~(\ref{eq:BZDataRefl}) follows from the isomorphisms
$$\Hom_{\Pi(\widehat Q)}\bigl(\widehat N(\gamma),T\bigr)\cong
\Hom_{\Pi(\widehat Q)}\bigl(\Sigma_i^*\widehat N(s_i\gamma),T\bigr)\cong
\Hom_{\Pi(\widehat Q)}\bigl(\widehat N(s_i\gamma),\Sigma_iT\bigr)$$
provided by Propositions~\ref{pr:NGammaFoncRefl} and~\ref{pr:ReflFuncAdj}.

Replacing $T$ by $\Sigma_i^*T$ in this equality~(\ref{eq:BZDataRefl})
yields
$$D_{\gamma}(\Sigma_i^*T)=D_{s_i\gamma}(\Sigma_i\Sigma_i^*T).$$
Now the equations (\ref{eq:DimHomNGamma}) say that
$$\dim\Hom_{\Pi(\widehat Q)}\bigl(\widehat N(s_i\gamma),S_i\bigr)=
-\langle\gamma,\alpha_i\rangle\quad\text{ and }\quad\Ext^1_{\Pi(\widehat
Q)}\bigl(\widehat N(s_i\gamma),S_i\bigr)=0.$$
Equation~(\ref{eq:InductRelPrep}) then follows by applying the
functor $\Hom_{\Pi(\widehat Q)}\bigl(\widehat N(s_i\gamma),?\bigr)$
to the short exact sequence
$$0\to\soc_iT\to T\to\Sigma_i\Sigma_i^*T\to0.$$

Finally, (\ref{eq:BZDataRadi}) is obtained by writing
(\ref{eq:InductRelPrep}) for the $\Pi(Q)$-module $\Sigma_iT$,
which has trivial $i$-socle.
\end{proof}

Let $a_{ij}=\langle\alpha_i,\alpha_j\rangle$ be the entries of the Cartan
matrix of the Dynkin diagram $\Gamma$.
\begin{proposition}
{\upshape\bfseries(Edge relations)}
\label{pr:EdgeRel}
Let $T$ be a $\Pi(Q)$-module. Then for each $i\in I$ and each $w\in W$,
one has
$$D_{-w\omega_i}(T)+D_{-ws_i\omega_i}(T)+\sum_{\substack{j\in I\\
j\neq i}}a_{ij}D_{-w\omega_j}(T)\geq0.$$
\end{proposition}
\begin{proof}
We fix $i\in I$ and call $L_w(T)$ the left hand side of the
desired relation. We want to show that $L_w(T)\geq0$ for all
$w\in W$. Since $L_w(T)=L_{ws_i}(T)$, we may restrict our
attention to the elements $w$ such that $\ell(ws_i)>\ell(w)$.

Take such a $w$, assume that it has positive length, and write
$w=s_kv$ with $k\in I$ and $\ell(v)=\ell(w)-1$.
Observing that $\ell(s_kv)>\ell(v)$ and $\ell(s_kvs_i)>
\ell(w)=\ell(v)+1\geq\ell(vs_i)$, a straightforward computation
based on equation (\ref{eq:InductRelPrep}) and on
$$\omega_i+s_i\omega_i+\sum_{\substack{j\in I\\j\neq i}}a_{ij}
\omega_j=0$$
shows that $L_w(T)=L_v(\Sigma_k^*T)$. By induction, we thus see
that if $s_{k_1}\cdots s_{k_r}$ is a reduced decomposition of $w$, then
$$L_w(T)=L_1(\Sigma_{k_r}^*\cdots\Sigma_{k_1}^*T).$$
Since $D_\gamma(T)=0$ for all antidominant weight $\gamma$, we finally
obtain
$$L_w(T)=D_{-s_i\omega_i}(\Sigma_{k_r}^*\cdots\Sigma_{k_1}^*T)\geq0.$$
\end{proof}

\begin{proposition}
\label{pr:CompRel}
Let $\gamma$ be a weight and $T$ be a $\Pi(Q)$-module. Then
$$D_\gamma(T)-D_{-\gamma}(T^*)=\langle\gamma,\dimvec T\rangle.$$
\end{proposition}
\begin{proof}
Let $i\in I$. By Proposition~\ref{pr:ProjCover},
$D_{\omega_i}(T)$ is equal to the Jordan-H\"older
multiplicity of $S_i$ in $T$, and therefore to
$\langle\omega_i,\dimvec T\rangle$. This fact implies the desired
equality in the case $\gamma=-\omega_i$. The case where $\gamma$ is
antidominant then follows from Theorem~\ref{th:CharNGamma}~\ref{it:CNGc}.

Fix an antidominant weight $\lambda$. For $w\in W$, call
$\mathscr E_w(T)$ the desired equality for $\gamma=w\lambda$.
Thus the first part of the proof establishes $\mathscr E_1(T)$
for any $T$.

Take $w\in W$ and $i\in I$ such that $\ell(s_iw)>\ell(w)$. The third
equality in Proposition~\ref{pr:NGammaFoncRefl} gives
$$\Sigma_i^*\widehat N(-w\lambda)\cong\widehat N(-s_iw\lambda).$$
By Proposition~\ref{pr:ReflFuncAdj} and Remark~\ref{rk:ReflFunc}
\ref{it:RFa}, it follows that
$$\Hom_{\Pi(\widehat Q)}\bigl(\widehat N(-s_iw\lambda),T^*\bigr)
\cong\Hom_{\Pi(\widehat Q)}\bigl(\widehat N(-w\lambda),\Sigma_iT^*\bigr)
\cong\Hom_{\Pi(\widehat Q)}\bigl(\widehat N(-w\lambda),(\Sigma_i^*T)^*\bigr),$$
whence
$$D_{-s_iw\lambda}(T^*)=D_{-w\lambda}\bigl((\Sigma_i^*T)^*\bigr).$$
A straightforward computation based on equation (\ref{eq:InductRelPrep})
and on
$$\dimvec T=\dimvec\Sigma_i\Sigma_i^*T+\dimvec\soc_iT=
s_i\bigl(\dimvec\Sigma_i^*T\bigr)+\dimvec\soc_iT$$
shows then that $\mathscr E_{s_iw}(T)$ is equivalent to
$\mathscr E_w(\Sigma_i^*T)$.

We conclude the proof by an immediate induction on $\ell(w)$.
\end{proof}

\begin{other}{Remark}
\label{rk:DGammaDom}
\begin{enumerate}
\item\label{it:DGDa}
It is worthwhile to record the starting point for the induction in
the proof of Proposition~\ref{pr:CompRel}: if $\gamma$ is dominant,
then
$$\dim\Hom_{\Pi(Q)}(N(\gamma),T)=D_\gamma(T)=\langle\gamma,\dimvec
T\rangle.$$
\item\label{it:DGDb}
Let us write a weight $\gamma$ as the difference $\gamma_0-\gamma_1$
of two dominant weights with disjoint ``support'': a fundamental
weight $\omega_i$ may appear in $\gamma_0$ or in $\gamma_1$, but not
in both. Let $P_0=N(\gamma_0)$ and $P_1=N(\gamma_1)$.
Proposition~\ref{pr:CompRel} then yields
$$\dim\Hom_{\Pi(Q)}(P_0,T)-\dim\Hom_{\Pi(Q)}(P_1,T)=
\langle\gamma,\dimvec T\rangle =D_\gamma(T)-D_{-\gamma}(T^*).$$
On the other hand, Proposition~\ref{pr:ProjCover} and
Theorem~\ref{th:CharNGamma}~\ref{it:CNGc} imply that $P_0$ and $P_1$
are projective modules, and examples seem to indicate that
$$\cdots\to P_1\to P_0\to N(\gamma)\to 0$$
is the beginning of the minimal projective resolution of $N(\gamma)$.
If this were true, the comparison of the formula just above with
Corollary~IV.4.3 in~\cite{Auslander-Reiten-Smalo95} would give
$$\dim\Hom_{\Pi(Q)}(T,D\,\mathrm{Tr}\,N(\gamma))=
D_{-\gamma}(T^*)=\dim\Hom_{\Pi(Q)}(T,N(-\gamma)^*).$$
One could then hope that $N(-\gamma)^*\cong D\,\mathrm{Tr}\,N(\gamma)$.
(Here, following the notation in~\cite{Auslander-Reiten-Smalo95}, 
$D $ and $ Tr $ are the $k$-duality and transpose functors 
respectively.
\end{enumerate}
\end{other}

\subsection{Pseudo-Weyl polytopes}
\label{ss:PWPoly}
Following Appendix~A in \cite{Kamnitzer05}, we recall the notion of
pseudo-Weyl polytope.

Let $V$ be a $\mathbb R$-vector space and let $V^*$ be its dual.
To a non-empty compact convex subset $P$ of $V$, we associate its
support function $\psi_P:V^*\to\mathbb R$: it maps a linear form
$\alpha\in V^*$ to the maximal value $\alpha$ takes on $P$.
Then $\psi_P$ is a sublinear function on $V^*$. One can recover
$P$ from the datum of $\psi_P$ by the Hahn-Banach theorem
$$P=\{v\in V\mid\forall\alpha\in V^*,\ \langle v,\alpha\rangle\leq
\psi_P(\alpha)\},$$
and the map $P\mapsto\psi_P$ is a bijection from the set of all
non-empty compact convex subsets of $V$ onto the set of all
sublinear functions on $V^*$ (see for instance Chapter~C in
\cite{Hiriart-Urruty-Lemarechal01}). If $P$ is a polytope, then
its support function is piecewise linear. More precisely, the
maximal regions of linearity of $\psi_P$ are exactly the maximal
cones of the dual fan of $P$: for each vertex $v$ of $P$, the
support function $\psi_P$ is linear on
$\{\alpha\in V^*\mid\psi_P(\alpha)=\langle v,\alpha\rangle\}$.

We now specialise to the case $V=Q\otimes_{\mathbb Z}\mathbb R$,
whence $V^*=P\otimes_{\mathbb Z}\mathbb R$. We say that a
polytope $P\subseteq V$ is a pseudo-Weyl polytope if its support
function is linear on each Weyl chamber. Then to each Weyl chamber
$C$ corresponds a vertex $v_C$ of $P$, defined by the condition
$$\forall\gamma\in C,\quad\psi_P(\gamma)=\langle\gamma,v_C\rangle.$$
One gets all the vertices of $P$ in this way (possibly with repetitions,
since the Weyl fan can be finer than the dual fan of $P$). For $w\in W$,
we denote by $\mu_w$ the vertex that corresponds to the chamber formed
by the elements $w\lambda$ with $\lambda$ antidominant. This indexed
collection $(\mu_w)_{w\in W}$ is called the vertex datum of $P$.

Let $\Gamma=\{w\omega_i\mid i\in I,\;w\in W\}$ denote the set of
chamber weights. Then each Weyl chamber is spanned by a subset of
$\Gamma$, which implies that the support function of a pseudo-Weyl
polytope is entirely characterized by its values on $\Gamma$.
A pseudo-Weyl polytope $P$ is thus characterized by the collection of
real numbers $\bigl(\psi_P(\gamma)\bigr)_{\gamma\in\Gamma}$;
concretely, these numbers describe the position of the facets of $P$:
$$P=\{v\in V\mid\forall\gamma\in\Gamma,\ \langle\gamma,v\rangle\leq
\psi_P(\gamma)\}.$$
This collection of real numbers is called the hyperplane datum of $P$.

Conversely, we can start from a collection of numbers
$(A_\gamma)_{\gamma\in\Gamma}$ and ask ourselves whether there
exists a sublinear function $\psi:V^*\to\mathbb R$, linear on
each Weyl chamber, such that $\psi(\gamma)=A_\gamma$ for each
$\gamma\in\Gamma$; in other words, whether
$(A_\gamma)_{\gamma\in\Gamma}$ is the hyperplane datum of a pseudo-Weyl
polytope. The answer is given in Lemma~A.5 in \cite{Kamnitzer05}:
a necessary and sufficient condition is that
$(A_\gamma)_{\gamma\in\Gamma}$ satisfy the edge inequalities
\begin{equation}
A_{-w\omega_i}+A_{-ws_i\omega_i}+\sum_{\substack{j\in I\\
j\neq i}}a_{ij}A_{-w\omega_j}\geq0
\label{eq:EdgeIneq}
\end{equation}
for all $i\in I$ and $w\in W$. When this condition is fulfilled,
the vertex datum $(\mu_w)_{w\in W}$ of the corresponding pseudo-Weyl
polytope is characterized by the set of equations
\begin{equation}
\langle-w\omega_i,\mu_w\rangle=A_{-w\omega_i},
\label{eq:RelVertHyp}
\end{equation}
and the left hand side of (\ref{eq:EdgeIneq}) is equal to the length $c$
of the edge between the vertices $\mu_w$ and $\mu_{ws_i}$, defined by
the equation $\mu_{ws_i}-\mu_w=cw\alpha_i$, by formula (8) in
\cite{Kamnitzer05}. From (\ref{eq:RelVertHyp}), we see that all the
weights $\mu_w$ belong to $Q$ if and only if all the numbers $A_\gamma$
are integers.

\subsection{GGMS strata}
\label{ss:GGMSStrata}
We keep the conventions of the previous section. For each
$\Pi(Q)$-module $T$, we have the collection of integers
$\big(D_\gamma(T)\big)_{\gamma\in\Gamma}$. Proposition
\ref{pr:EdgeRel} shows that this collection satisfies the edge
inequalities. It is thus the hyperplane datum of a pseudo-Weyl
polytope, namely
$$\Pol(T)=\{v\in V\mid\forall\gamma\in\Gamma,\ \langle\gamma,v
\rangle\leq D_\gamma(T)\}.$$
Thus $D_\gamma(T)$ is the value $\psi_{\Pol(T)}(\gamma)$ of the
support function, for each $\gamma\in\Gamma$, and indeed for each
$\gamma\in P$, thanks to Theorem~\ref{th:CharNGamma}~\ref{it:CNGc}.
The vertex datum $\bigl(\mu_w(T)\bigr)_{w\in W}$ of $\Pol(T)$ can
be obtained from the numbers $D_\gamma(T)$ by the formula
(\ref{eq:RelVertHyp}), with $D_\gamma(T)$ instead of $A_\gamma$. Using
Remark~\ref{rk:DGammaDom}~\ref{it:DGDa}, we note that $\mu_1(T)=0$
and $\mu_{w_0}(T)=\dimvec T$.

Conversely, for each dimension-vector $\nu$ and each pseudo-Weyl
polytope $P$ whose vertex datum satisfies $\mu_1=0$ and $\mu_{w_0}=\nu$,
we may consider the set of all points in $\Lambda(\nu)$ whose
polytope is $P$. We obtain in this way a stratification of
$\Lambda(\nu)$ whose strata are labelled by pseudo-Weyl
polytopes of weight $\nu$. In analogy with the situation for the
affine Grassmannian (see section 2.4 of \cite{Kamnitzer05}), we
call these the GGMS strata of $\Lambda(\nu)$. Later in this paper,
we will see that the strata of biggest dimension (whose closures
are the components) are those labelled by MV polytopes ---
exactly as in the case of the affine Grassmannian. On the other
hand, we do not know how to analyse the other GGMS strata and
compare them with the corresponding GGMS strata in the affine
Grassmannian.

To conclude, let us emphasize the meaning of Proposition~\ref{pr:CompRel}:
for any $\Pi(Q)$-module $T$, the pseudo-Weyl polytope $\Pol(T^*)$ is
the image of $\Pol(T)$ under the involution $x\mapsto\dimvec T-x$.

\section{Crystal structure}
\label{se:CrysStruct}
The combinatorics of representations of complex semi-simple Lie algebras
is in part controlled by Kashiwara crystals; see \cite{Kashiwara95}
for a survey of this theory. The crystals $B(\lambda)$ of the
finite-dimensional representations of $\mathfrak g$ are ``contained''
in a big crystal $B(-\infty)$, which is the crystal of the positive
part $U(\mathfrak n)$ of $U(\mathfrak g)$. In \cite{Kashiwara-Saito97},
Kashiwara and Saito gave a characterization of $B(-\infty)$ as the
unique highest weight crystal having an involution with certain
specific properties.

Let us denote by $\Irr X$ the set of irreducible components of an
algebraic variety $X$. The set $\mathbf B=\bigsqcup_{\nu\in Q_+}
\Irr\Lambda(\nu)$ can be endowed with the structure of a crystal
with an involution, which turns out to be isomorphic to $B(-\infty)$.
Thus to an element $b\in B(-\infty)$, of weight $\nu$, corresponds
$\Lambda_b\in\Irr\Lambda(\nu)$. Our aim in this section is to show that
the reflection functor $\Sigma_i$ has a crystal counterpart, namely
the map denoted by $S_i$ in \cite{Kashiwara-Saito97}. (We will adopt
this notation $S_i$, though it brings some confusion with our
previous notation for the simple quiver representations.)

\subsection{The crystal structure}
\label{ss:CrysStruct}
We begin with recalling the crystal structure on $\mathbf B$, first
defined by Lusztig in \cite{Lusztig90b}, section~8.

Let $\nu\in Q_+$. As in section~\ref{ss:PreprojAlg}, we set
$M_i=K^{\nu_i}$ for each $i\in I$. Any point $(M_a)\in\Lambda(\nu)$
then defines a $\Pi(Q)$-module with dimension-vector $\nu$. An
element $Z\in\Irr\Lambda(\nu)$ is said to have weight $\wt(Z)=\nu$.

The usual isomorphism between the vector spaces $M_i$ and their
duals allows to view the involution $*$ on $\Pi(Q)$-modules as an
operation on $\Lambda(\nu)$ (it just amounts to the transposition
and subsequent relabelling of the matrices $M_a$). This operation
induces an involution on $\Irr\Lambda(\nu)$. We thus obtain a weight
preserving involution on $\mathbf B$, which we denote again by $*$.

For $i\in I$ and $c\in\mathbb N$, let
$$\Lambda(\nu)_{i,c}=\bigl\{(M_a)\bigm|\dim\hd_iM=c\bigr\}.$$
These sets form a partition of $\Lambda(\nu)$ into locally closed
subsets. For any irreducible component $Z$ of $\Lambda(\nu)$, there
is thus one value of $c$ such that $Z\cap\Lambda(\nu)_{i,c}$ is
open and dense in $Z$. We then write $\varphi_i(Z)=c$ and
$\varepsilon_i(Z)=c-\langle\alpha_i,\nu\rangle$.

In this context, we note that $\Lambda(\nu)$ and each subset
$\Lambda(\nu)_{i,c}$ have pure dimension equal to
$\dim G(\nu)-(\nu,\nu)/2$; see Theorem~8.7 in~\cite{Lusztig90b}.
Therefore the map $Z\mapsto Z\cap\Lambda(\nu)_{i,c}$ gives a bijection
$$\{Z\in\Irr\Lambda(\nu)\mid\varphi_i(Z)=c\}\to\Irr\Lambda(\nu)_{i,c}.$$

Let again $i\in I$. Given $c\in\mathbb N$, let us denote
by $\Omega(\nu,i,c)$ the set of triples $((M_a),(N_a),g)$
such that $(M_a)\in\Lambda(\nu)_{i,0}$,
$(N_a)\in\Lambda(\nu+c\alpha_i)_{i,c}$ and $g:M\to N$ is an injective
morphism of $\Pi(Q)$-modules. We can then form the diagram
\begin{equation}
\Lambda(\nu)_{i,0}\xleftarrow p\Omega(\nu,i,c)\xrightarrow q
\Lambda(\nu+c\alpha_i)_{i,c},
\label{eq:DiagCrysStruct}
\end{equation}
where $p$ and $q$ are the first and second projection.
Then $p$ is a locally trivial fibration, with a smooth and connected
fiber, and $q$ is a principal $G(\nu)$-bundle. Thus $p$ and $q$
define bijections
$$\Irr\Lambda(\nu)_{i,0}\longleftrightarrow\Irr\Omega(\nu,i,c)
\longleftrightarrow\Irr\Lambda(\nu+c\alpha_i)_{i,c}.$$
We therefore obtain mutually inverse bijections
$$\xymatrix@C=3.7em{\{Z\in\Irr\Lambda(\nu)\mid\varphi_i(Z)=0\}
\ar@<.5ex>[r]^(.45){\tilde e_i^c}&
\{Z\in\Irr\Lambda(\nu+c\alpha_i)\mid\varphi_i(Z)=c\},
\ar@<.5ex>[l]^(.55){\tilde f_i^{\max{}}}}$$
whose data is equivalent to the data of the usual structure
maps for a crystal
$$\xymatrix@C=3em{\{Z\in\Irr\Lambda(\nu+c\alpha_i)\mid\varphi_i(Z)=c\}
\ar@<.5ex>[r]^(.45){\tilde e_i}&
\{Z\in\Irr\Lambda(\nu+(c+1)\alpha_i)\mid\varphi_i(Z)=c+1\}.
\ar@<.5ex>[l]^(.55){\tilde f_i}}$$

It is known that the maps $\wt$, $\varepsilon_i$, $\varphi_i$,
$\tilde e_i$ and $\tilde f_i$ endow $\mathbf B$ with the structure
of a crystal. Using the involution $*$, Kashiwara and Saito proved
the existence of an isomorphism $b\mapsto\Lambda_b$
from $B(-\infty)$ onto $\mathbf B$ (Theorem~5.3.2 in
\cite{Kashiwara-Saito97}). This isomorphism is unique since $B(-\infty) $ has no automorphism. As usual, one writes
$\tilde e_i^*=*\tilde e_i*$ and $\tilde f_i^*=*\tilde f_i*$.

For $\nu\in Q_+$, $i\in I$ and $c\in\mathbb N$, let us set
$\Lambda(\nu)_{i,c}^\times=\Lambda(\nu)_{i,c}\cap(\Lambda(\nu)_{i,0})^*$.
Assuming moreover that $0\leq c\leq-\langle\alpha_i,\nu\rangle$, let us
denote by $\Omega(\nu,i,c)^\times$ the set of triples $((M_a),(N_a),g)$
such that $(M_a)\in\Lambda(\nu)_{i,0}^\times$,
$(N_a)\in\Lambda(\nu+c\alpha_i)_{i,c}^\times$ and $g:M\to N$ is an
injective morphism of $\Pi(Q)$-modules. By restriction, the diagram
(\ref{eq:DiagCrysStruct}) yields
\begin{equation}
\Lambda(\nu)_{i,0}^\times\xleftarrow{p^\times}\Omega(\nu,i,c)^\times
\xrightarrow{q^\times}\Lambda(\nu+c\alpha_i)_{i,c}^\times.
\label{eq:DiagCrysStructTimes}
\end{equation}

\begin{lemma}
\label{le:RestDiag}
In the context above, $p^\times$ is a locally trivial fibration, with
a smooth and connected fiber, and $q^\times$ is a principal
$G(\nu)$-bundle.
\end{lemma}
\begin{proof}
Let us set $M_i=K^{\nu_i}$, $N_i=K^{\nu_i+c}$, and $M_j=N_j=K^{\nu_j}$
for $j\neq i$.

Let us first study the fibers of $p^\times$. We thus fix a point
$(M_a)\in\Lambda(\nu)_{i,0}^\times$ and look at the set of
pairs $((N_a),g)$, where $(N_a)\in\Lambda(\nu+c\alpha_i)_{i,c}^\times$
and $g:M\to N$ is an injective morphism. The isomorphisms of vector
spaces $g_j:M_j\to N_j$ with $j\neq i$ can be freely chosen, and their
datum determines the maps $N_a$ for all $a\in H$ such that
$i\notin\{s(a),t(a)\}$. It thus remains to study the local situation
around $i$.

We adopt the notations set up in section~\ref{ss:DefReflFunc}.
The situation can be summarised by a commutative diagram
$$\xymatrix@C=3.7em{\widetilde M_i\ar[d]_{\widetilde g_i}\ar@{->>}
[r]^{M_{\iin(i)}}&M_i\ar@{^{(}->}^{g_i}[d]\ar@{^{(}->}[r]^{M_{\out(i)}}
&\widetilde M_i\ar[d]^{\widetilde g_i}\\\widetilde N_i
\ar[r]_{N_{\iin(i)}}&N_i\ar@{^{(}->}[r]_{N_{\out(i)}}&\widetilde N_i,}$$
where $\widetilde g_i$ is an isomorphism of vector spaces. The
injectivity of $M_{\out(i)}$ and $N_{\out(i)}$ follows from the
fact that $(M_a)\in(\Lambda(\nu)_{i,0})^*$ and $(N_a)\in(\Lambda
(\nu+c\alpha_i)_{i,0})^*$.

The image of the map $u=M_{\out(i)}M_{\iin(i)}$ has the same dimension
as $M_i$, namely $\nu_i$. The image of $(\widetilde g_i)^{-1}N_{\out(i)}$
is a subspace $V$ of $\widetilde M_i$ with dimension $\nu_i+c$ and such
that $\im u\subseteq V\subseteq\ker u$. In other words, $V/\im u$ can
be chosen freely in the Grassmannian of $c$-dimensional subspaces of
$\ker u/\im u$ (note that $\dim\ker u/\im u=-\langle\alpha_i,\nu\rangle
\geq c$). The choice of such a $V$ and of an isomorphism of vector
spaces $\overline N_{\out(i)}:N_i\to\widetilde g_i(V)$ determines the
remaining data $g_i$, $N_{\iin(i)}$ and $N_{\out(i)}$. We conclude that
the fiber of $p$ above $(M_a)$ has the structure of a principal
$G(\nu+c\alpha_i)$-bundle over a Grassmannian. This Grassmannian depends
smoothly on $(M_a)$, so we can conclude that $p^\times$ is a locally
trivial fibration, with a smooth and connected fiber.

The proof of the assertion about $q^\times$ is similar, but simpler.
\end{proof}

We can now write a commutative diagram
$$\xymatrix{\Irr\Lambda(\nu)_{i,0}^\times\ar@{^{(}->}[d]\ar@{<->}[r]&
\Irr\Omega(\nu,i,c)^\times\ar@{^{(}->}[d]\ar@{<->}[r]&
\Irr\Lambda(\nu+c\alpha_i)_{i,c}^\times\ar@{^{(}->}[d]\\
\Irr\Lambda(\nu)_{i,0}\ar@{<->}[r]&
\Irr\Omega(\nu,i,c)\ar@{<->}[r]&
\Irr\Lambda(\nu+c\alpha_i)_{i,c},}$$
where the horizontal arrows are the bijections defined by the maps
in the diagrams (\ref{eq:DiagCrysStructTimes}) and
(\ref{eq:DiagCrysStruct}) and the vertical arrows are given by the
open inclusions
$$\Lambda(\nu)_{i,0}^\times\subseteq\Lambda(\nu)_{i,0},\quad
\Omega(\nu,i,c)^\times\subseteq\Omega(\nu,i,c)\quad\text{and}\quad
\Lambda(\nu+c\alpha_i)_{i,c}^\times\subseteq\Lambda(\nu+c\alpha_i)_{i,c}.$$
In particular, we see that if $0\leq c\leq-\langle\alpha_i,\nu\rangle$,
then $\tilde e_i^c$ and $\tilde f_i^{\max{}}$ restrict to bijections
$$\xymatrix@C=2.2em{\{Z\in\Irr\Lambda(\nu)\mid\varphi_i(Z)=
\varphi_i(Z^*)=0\}\ar@<.5ex>[r]^(.435){\tilde e_i^c}&
\{Z\in\Irr\Lambda(\nu+c\alpha_i)\mid\varphi_i(Z)=c\text{ and
}\varphi_i(Z^*)=0\}.\ar@<.5ex>[l]^(.565){\tilde f_i^{\max{}}}}$$

\begin{other*}{Remark}
Taking into account the crystal isomorphism $\mathbf B\cong B(-\infty)$,
we conclude that for any $b\in B(-\infty)$ and $c\in\mathbb N$,
$$\bigl(\;\varphi_i(b)=\varphi_i(b^*)=0\quad\text{and}\quad
0\leq c\leq-\langle\alpha_i,\wt(b)\rangle\;\bigr)\ \Longrightarrow\
\varphi_i(\,(\tilde e_i^cb)^*\,)=0.$$
This result is indeed contained in Proposition~5.3.1~(1) in
\cite{Kashiwara-Saito97}, which says that
$$\varphi_i(b)=\max\bigl(\,\varphi_i\bigl((\tilde f_i^*)^{\max{}}b\bigr),
\;\langle\wt(b),\alpha_i\rangle-\varphi_i(b^*)\,\bigr)$$
for each $b\in B(-\infty)$. A particular case of the latter formula
is that for each $b\in B(-\infty)$ that satisfies $\varphi_i(b)=0$, the
number $\varepsilon_i(b^*)=\varphi_i(b^*)-\langle\alpha_i,\wt(b)\rangle$
is non-negative.
\end{other*}

\subsection{Reflection functors and crystal operations}
\label{ss:ReflFuncCrysOp}
Let us fix $i\in I$ for this whole section. In Corollary~3.4.8 in
\cite{Saito94} (see also section~8.2 in \cite{Kashiwara-Saito97}),
Saito defines a bijection
$$S_i:\{b\in B(-\infty)\mid\varphi_i(b)=0\}\to
\{b\in B(-\infty)\mid\varphi_i(b^*)=0\}$$
by the rule $S_i(b)=\tilde e_i^{}{}^{\!\!\varepsilon_i(b^*)}(\tilde
f_i^*)^{\max{}}b$. This operation plays the role of the simple
reflection $s_i$, but at the crystal level. One can indeed check
that
$$\wt(S_ib)=\wt(b)+(\varepsilon_i(b^*)-\varphi_i(b^*))\,\alpha_i=
\wt(b)-\langle\alpha_i,\wt(b^*)\rangle\,\alpha_i=s_i(\wt(b)).$$
A deeper property is that $S_i$ implements at the crystal level the
action of Lusztig's braid group automorphism $T_i=T'_{i,-1}$ on the
canonical basis; see \cite{Lusztig96} and Proposition~3.4.7 in
\cite{Saito94}. We now relate this operation $S_i$ to the reflection
functor $\Sigma_i$.

To do that, let us fix $\nu\in Q_+$ and let us denote by $\Theta(\nu,i)$
the set of all triples $((M_a),(N_a),h)$ such that
$(M_a)\in\Lambda(\nu)_{i,0}$, $(N_a)\in(\Lambda(s_i\nu)_{i,0})^*$
and $h:N\to\Sigma_iM$ is an isomorphism. We can then form the diagram
\begin{equation}
\Lambda(\nu)_{i,0}\xleftarrow r\Theta(\nu,i)\xrightarrow s
(\Lambda(s_i\nu)_{i,0})^*,
\label{eq:DiagReflFuncCrysOp}
\end{equation}
where $r$ and $s$ are the first and second projection.
From the definition, it is immediate that $r$ is a principal
$G(s_i\nu)$-bundle. On the other hand $s$ is a principal
$G(\nu)$-bundle, because any $(N_a)\in(\Lambda(s_i\nu)_{i,0})^*$ has
trivial $i$-socle and hence is in the essential image of $\Sigma_i$
by the remarks at the end of section~\ref{ss:DefReflFunc}. Thus $r$
and $s$ are locally trivial fibrations with a smooth and connected
fiber; therefore they define bijections
$$\Irr\Lambda(\nu)_{i,0}\longleftrightarrow\Irr\Theta(\nu,i)
\longleftrightarrow\Irr\bigl((\Lambda(s_i\nu)_{i,0})^*\bigr).$$

\begin{theorem}
\label{th:ReflFuncCrysOp}
Let $b\in B(-\infty)$ of weight $\nu$ and such that $\varphi_i(b)=0$.
In the notation above,
$$r^{-1}\bigl(\Lambda_b\cap\Lambda(\nu)_{i,0}\bigr)=
s^{-1}\bigl(\Lambda_{S_ib}\cap(\Lambda(s_i\nu)_{i,0})^*\bigr).$$
\end{theorem}
\begin{proof}
Let $\mu\in Q_+$ and let $c,d\in\mathbb N$ be such that
$c+d=-\langle\alpha_i,\mu\rangle$. Set $\nu=\mu+c\alpha_i$; thus
$s_i\nu=\mu+d\alpha_i$. A direct calculation shows that for
any $\Pi(Q)$-modules $M$ and $N$ such that
$$\dimvec M=\nu,\quad\hd_iM=0,\quad\dimvec N=s_i\nu,\quad
\soc_iN=0,\quad\Sigma_iM\cong N,$$
one has
$$\dim\soc_iM=c\ \Longleftrightarrow\ \dim\hd_iN=d.$$
The diagram (\ref{eq:DiagReflFuncCrysOp}) thus restrict to
\begin{equation}
\Lambda(\nu)_{i,0}\cap(\Lambda(\nu)_{i,c})^*
\leftarrow\Theta(\nu,i)_{c,d}\rightarrow
(\Lambda(s_i\nu)_{i,0})^*\cap\Lambda(s_i\nu)_{i,d},
\label{eq:DiagReflFuncCrysOpTimes}
\end{equation}
for a suitable locally closed subset
$\Theta(\nu,i)_{c,d}\subseteq\Theta(\nu,i)$.

Now let $\Xi$ be the set of all tuples
$((L_a),(M_a),(N_a),f,g,h)$ such that
$$(L_a)\in\Lambda(\mu)_{i,0}^\times,\quad
(M_a)\in(\Lambda(\nu)_{i,c}^\times)^*,\quad
(N_a)\in\Lambda(s_i\nu)_{i,d}^\times,$$
and that $f:M\twoheadrightarrow L$,\quad$g:L\hookrightarrow N$
and $h:N\stackrel\simeq\longrightarrow\Sigma_iM$
are morphisms of $\Pi(Q)$-modules,
subject to the condition that $hgf:M\to\Sigma_iM$ is the canonical
map from Remark~\ref{rk:ReflFunc}~\ref{it:RFb}.

We inscribe (\ref{eq:DiagReflFuncCrysOpTimes}) as the right column in
$$\xymatrix@C=4em{(\Lambda(\mu)_{i,0}^\times)^*\ar@{=}[dd]&
(\Omega(\mu,i,c)^\times)^*
\ar[l]_(.53){*\;p^\times*}\ar[r]^(.53){*\;q^\times*}&
(\Lambda(\nu)_{i,c}^\times)^*\\
&\Xi\ar[u]\ar[r]\ar[d]&\Theta(\nu,i)_{c,d}\ar[u]_r\ar[d]^s\\
\Lambda(\mu)_{i,0}^\times&
\Omega(\mu,i,d)^\times\ar[l]^{p^\times}\ar[r]_(.53){q^\times}&
\Lambda(s_i\nu)_{i,d}^\times.}$$
In this commutative diagram, the arrows are the relevant projections
and the top and bottom rows are diagrams of the form
(\ref{eq:DiagCrysStructTimes}). Routine arguments show that the three
arrows starting from $\Xi$ are principal bundles, with structural
groups $G(s_i\nu)$, $G(\mu)$ and $G(\nu)$. Therefore all the maps on
this diagram are locally trivial fibrations with a smooth and
connected fiber. They thus induce bijections between the sets of
irreducible components.

Now take $b\in B(-\infty)$ as in the statement of the theorem. Let
$c=\varphi_i(b^*)$, $d=\varepsilon_i(b^*)$ and
$\mu=\wt\bigl((\tilde f_i^*)^{\max{}}b\bigr)$;
then the weight of $b$ is $\nu=\mu+c\alpha_i$. The remark at the
end of section~\ref{ss:CrysStruct} says that
$d=c-\langle\alpha_i,\nu\rangle$ is non-negative and that
$\varphi_i\bigl((\tilde f_i^*)^{\max{}}b\bigr)=0$, which implies
$$d=\varphi_i(S_ib)\quad\text{and}\quad
c+d=2c-\langle\alpha_i,\mu+c\alpha_i\rangle=-\langle\alpha_i,\mu\rangle.$$
This remark also says that $\varphi_i\bigl((S_ib)^*\bigr)=0$.
The commutativity of our big diagram implies that
$\Lambda_b\cap(\Lambda(\nu)_{i,c}^\times)^*$ and
$\Lambda_{S_ib}\cap\Lambda(s_i\nu)_{i,d}^\times$ correspond
in the bijection defined by (\ref{eq:DiagReflFuncCrysOpTimes}).
We conclude by taking closures in $\Lambda(\nu)_{i,0}$ and
$(\Lambda(s_i\nu)_{i,0})^*$.
\end{proof}

\subsection{Reformulations}
\label{ss:Reform}
In this subsection, we rewrite in a more direct way the
constructions presented in sections~\ref{ss:CrysStruct} and
\ref{ss:ReflFuncCrysOp}.

Given a dimension vector $\nu\in Q_+$, any point $(M_a)\in\Lambda(\nu)$
can be viewed as a $\Pi(Q)$-module with dimension-vector $\nu$. Conversely,
if $M$ is a $\Pi(Q)$-module with dimension-vector $\nu$, then its
isomorphism class can be viewed as a $G(\nu)$-orbit in $\Lambda(\nu)$;
this orbit will be denoted by $[M]$.

We now fix $\nu\in Q_+$ and $i\in I$. Our first reformulation deals
with the operation $\tilde f_i^{\max{}}$.

\begin{proposition}
\label{pr:ReformFimax}
Let $Z$ be an irreducible component of $\Lambda(\nu)$ and set
$c=\varphi_i(Z)$ and $\nu'=\nu-c\alpha_i$. Let $V$ be a dense, open,
and $G(\nu')$-invariant subset of $\tilde f_i^{\max{}}Z$. Then
$\{M\in Z\mid[\Sigma_i^*\Sigma_iM]\subseteq V\}$
contains a dense open subset in $Z$.
\end{proposition}
\begin{proof}
Let $c\in\mathbb N$ and set $\nu'=\nu-c\alpha_i$. The proof of
Proposition~\ref{pr:ReflFuncAdj}~\ref{it:RFAb} tells us that if
$M$ is a $\Pi(Q)$-module with dimension-vector $\nu$ with an $i$-head
of dimension $c$, then $\Sigma_i^*\Sigma_iM$ identifies canonically
with the unique submodule of $M$ with dimension-vector $\nu'$. With
the notation of section~\ref{ss:CrysStruct}, this means that
$\Omega(\nu',i,c)$ is the set of triples $((N_a),(M_a),g)$ such that
$(N_a)\in\Lambda(\nu')_{i,0}$, $(M_a)\in\Lambda(\nu)_{i,c}$ and
$g:N\to\Sigma_i^*\Sigma_iM$ is an isomorphism of $\Pi(Q)$-modules.
Looking at the diagram (\ref{eq:DiagCrysStruct}), we thus see that
for each $(M_a)\in\Lambda(\nu)_{i,c}$, the set $p(q^{-1}(M))$ is
the orbit $[\Sigma_i^*\Sigma_iM]$.

Now let $Z$ be an irreducible component of $\Lambda(\nu)$. Put
$c=\varphi_i(Z)$ in the discussion just above; thus $U=Z\cap\Lambda(\nu)_{i,c}$
is a dense open subset of $Z$. By definition of the operation
$\tilde f_i^{\max{}}$, the diagram~(\ref{eq:DiagCrysStruct})
restricts to
$$U'\xleftarrow p\Omega'\xrightarrow qU,$$
where $\Omega'$ is an irreducible component of $\Omega(\nu,i,c)$
and $U'=\bigl(\tilde f_i^{\max{}}Z\bigr)\cap\Lambda(\nu')_{i,0}$ is
dense open in $\tilde f_i^{\max{}}Z$. From the fact that $U'\cap V$
is open dense in $U'$, it follows that $q(p^{-1}(U'\cap V))$ is
open dense in $U$.

To conclude the proof, it now suffices to observe that
\begin{align*}
\{M\in Z\mid[\Sigma_i^*\Sigma_iM]\subseteq V\}
&=\{M\in Z\mid[\Sigma_i^*\Sigma_iM]\text{ meets }V\}\\
&\supseteq\{M\in U\mid p(q^{-1}(M))\text{ meets }U'\cap V\}\\
&=q(p^{-1}(U'\cap V));
\end{align*}
in this computation, the first equality comes from the
$G(\nu')$-invariance of $V$.
\end{proof}

Our second reformulation concerns the operation $S_i$.

\begin{proposition}
\label{pr:ReflFuncGene}
Let $b\in B(-\infty)$ and set
$$\nu=\wt(b),\quad c=\varphi_i(b),\quad
U=\Lambda_b\cap\Lambda(\nu)_{i,c},\quad
b'=S_i\bigl(\tilde f_i^{\max{}}b\bigr)\quad\text{and}\quad
\nu'=\wt(b').$$
Let $V$ be an open, dense and $G(\nu')$-invariant subset of
$\Lambda_{b'}$. Then $\{(M_a)\in U\mid[\Sigma_iM]\subseteq V\}$
contains a dense open subset of $U$.
\end{proposition}
\begin{proof}
The proof of the particular case $c=0$ is completely analogous to
the proof of Proposition~\ref{pr:ReformFimax}. Indeed we first
note that $b'=S_ib$, hence $\nu'=s_i\nu$, and then we let the
diagram (\ref{eq:DiagReflFuncCrysOp}) play the role of the diagram
(\ref{eq:DiagCrysStruct}) in our previous discussion, using
Theorem~\ref{th:ReflFuncCrysOp} to identify the bijection between
the sets of irreducible components.

The general case is a combination of this particular case $c=0$ with
the situation studied in Proposition~\ref{pr:ReformFimax}. This is
obvious on the crystal side; on the side of $\Pi(Q)$-modules, one uses
the equality $\Sigma_iM\cong\Sigma_i(\Sigma_i^*\Sigma_iM)$ noticed in
section~\ref{ss:DefReflFunc}.
\end{proof}

\section{MV polytopes}
\label{se:MVPoly}
In section~\ref{ss:GGMSStrata}, to each $\Pi(Q)$-module $T$, we
associated the pseudo-Weyl polytope $\Pol(T)$ with hyperplane datum
$\bigl(D_\gamma(T)\bigr)_{\gamma\in\Gamma}$. MV polytopes are
special pseudo-Weyl polytopes, defined by certain rank $2$ conditions.
The work of Anderson~\cite{Anderson03} and the second author
\cite{Kamnitzer05} establishes a link between these MV polytopes and
the combinatorics of representations. In particular, MV polytopes
form a model for the crystal $B(-\infty)$.

Our goal in this section is to show that if $T$ is a general point
in a component $Z\subseteq\Lambda(\nu)$, then the polytope $\Pol(T)$
constructed in section~\ref{ss:GGMSStrata} is an MV polytope.
Moreover, we show that this construction provides an isomorphism of
crystals between $\mathbf B$ and the set of MV polytopes.

\subsection{Tropical Pl\"ucker relations}
\label{ss:TropPlucRel}
A hexagon is a triple $(w,i,j)\in W\times I^2$ such that
$a_{ij}=-1$, $\ell(ws_i)>\ell(w)$ and $\ell(ws_j)>\ell(w)$.
(This terminology comes from the fact that such a triple corresponds to
an hexagonal $2$-face of any pseudo-Weyl polytope.)

We say that a collection of real numbers $(A_\gamma)_{\gamma\in\Gamma}$
satisfies the tropical Pl\"ucker relation at the hexagon $(w,i,j)$ if
$$A_{-ws_i\omega_i}+A_{-w s_j\omega_j}=\max(A_{-w\omega_i}+
A_{-ws_i s_j\omega_j},A_{-ws_js_i\omega_i}+A_{-w\omega_j}).$$
(The minus signs here are not essential, for they can be removed
by replacing $w$ by $ww_0$. They have been inserted to simplify the
statement of Lemma~\ref{le:TropPluRelBase}.)

We say that a collection $(A_\gamma)_{\gamma\in\Gamma}$ is a BZ datum
if the following three conditions hold:
\vspace{-7pt}
\begin{description}
\item[(BZ1)]
Each $A_\gamma$ is an integer and each $A_{-\omega_i}$ is zero.
\item[(BZ2)]
All edge inequalities (\ref{eq:EdgeIneq}) are satisfied.
\item[(BZ3)]
All possible tropical Pl\"ucker relations are satisfied.
\end{description}
\vspace{-7pt}
We say that a pseudo-Weyl polytope is an MV polytope if its hyperplane
datum is a BZ datum. Thus there is a canonical bijection between the
set of MV polytopes and the set of BZ data.

For any $\Pi(Q)$-module $T$, the collection
$\bigl(D_\gamma(T)\bigr)_{\gamma\in\Gamma}$ obviously satisfies (BZ1).
Proposition~\ref{pr:EdgeRel} says that it also satisfies (BZ2); this
fact was indeed used in section~\ref{ss:GGMSStrata} to associate to
$T$ its pseudo-Weyl polytope $\Pol(T)$. We will now prove that when
$T$ is the general point in a component of a variety $\Lambda(\nu)$,
then $\bigl(D_\gamma(T)\bigr)_{\gamma\in\Gamma}$ satisfies (BZ3);
in other words, $\Pol(T)$ is an MV polytope.

More precisely, note that $D_\gamma$ is a constructible function on
$\Lambda(\nu)$. Hence any irreducible component $Z$ of $\Lambda(\nu)$
contains a dense open subset on which $D_\gamma$ takes a constant
value: we denote this value by $D_\gamma(Z)$. There is thus an open
and dense subset $\Omega_Z\subseteq Z$ such that each point $T\in\Omega_Z$
satisfies $D_\gamma(T)=D_\gamma(Z)$ for each $\gamma\in\Gamma$
(and hence for each $\gamma\in P$); and one may even demand that
$\Omega_Z$ be $G(\nu)$-invariant, for each function $D_\gamma$ is
$G(\nu)$-invariant.

\begin{theorem}
\label{th:TropPluRel}
Let $\nu\in\mathbb N^I$ be a dimension-vector. Then for any irreducible
component $Z$ of $\Lambda(\nu)$, the collection $\bigl(D_\gamma(Z)
\bigr)_{\gamma\in\Gamma}$ satisfies the tropical Pl\"ucker relations.
\end{theorem}

We have to prove a relation at each hexagon $(w,i,j)$. We will proceed
by induction on $\ell(w)$. We begin by checking the base case $w=1$.

\begin{lemma}
\label{le:TropPluRelBase}
Let $(i,j)$ be a pair of vertices in the Dynkin diagram, connected
to each other. Let $\nu$ be a dimension-vector and $Z$ be an
irreducible component of $\Lambda(\nu)$. Then
$\bigl(D_\gamma(Z)\bigr)$ satisfies the tropical Pl\"ucker
relation at $(1,i,j)$.
\end{lemma}
\begin{proof}
Let $i$, $j$ and $Z$ be as in the statement of the lemma.
We first note that $N(-s_i\omega_i)\cong S_i$ and $N(-s_j\omega_j)\cong
S_j$. With the diagrammatical convention of section~\ref{ss:PreprojAlg},
the $\Pi(Q)$-modules $T_i=N(-s_js_i\omega_i)\cong\Sigma_jS_i$ and
$T_j=N(-s_is_j\omega_j)\cong\Sigma_iS_j$ are represented by
$$T_i:\ \raisebox{-1.6ex}{\xymatrix@ur@=1.3em{i&j\ar[l]}}\qquad\text{and}
\qquad T_j:\ \raisebox{2.7ex}{\xymatrix@dr@=1.3em{i\ar[r]&j.}}$$

We want to show that for a general point $M$ in $Z$,
$$\dim\Hom(S_i,M)+\dim\Hom(S_j,M)=
\max\bigl(\dim\Hom(T_i,M),\dim\Hom(T_j,M)\bigr).$$

So let us take a point $M$ in $\Omega_Z$, and let us adopt the same
local representation for $M$ as the one used in the proof of
Proposition~\ref{pr:BraidRel}: we abbreviate the part
\vspace{-4ex}
$$\xymatrix@C=5em{
\raisebox{-7ex}{$\displaystyle\bigoplus_{\substack{a\in
H\\t(a)=i\\a\neq c^*}}M_{s(a)}$}
\ar@<.5ex>[r]^(.6){(\varepsilon(a)M_a)}&
M_i\ar@<.5ex>[l]^(.4){(M_{a^*})}\ar@<.5ex>[r]^{M_c}&
M_j\ar@<.5ex>[l]^{M_{c^*}}\ar@<.5ex>[r]^(.4){(M_{b^*})}&
\raisebox{-7ex}{$\displaystyle\bigoplus_{\substack{b\in
H\\t(b)=j\\b\neq c}}M_{s(b)}$}
\ar@<.5ex>[l]^(.6){(\varepsilon(b)M_b)}}$$
of the datum of $M$ by the notation
$$\xymatrix{R\ar@<.5ex>[r]^k&V\ar@<.5ex>[l]^l\ar@<.5ex>[r]^f&
W\ar@<.5ex>[l]^g\ar@<.5ex>[r]^m&S.\ar@<.5ex>[l]^n}$$
The preprojective relations at $i$ and $j$ then read
$$kl+\varepsilon(c^*)gf=0\quad\text{ and }\quad
nm+\varepsilon(c)fg=0.$$

We are interested in the dimension of the spaces
\begin{align*}
\Hom(S_i,M)&\cong\soc_iM=\ker f\cap\ker l,\\
\Hom(S_j,M)&\cong\soc_jM=\ker g\cap\ker m,\\
\Hom(T_i,M)&\cong\ker m\cap g^{-1}(\soc_iM)
=\ker m\cap g^{-1}(\ker f\cap\ker l),\\
\Hom(T_j,M)&\cong\ker l\cap f^{-1}(\soc_jM)
=\ker l\cap f^{-1}(\ker g\cap\ker m).
\end{align*}

By the preprojective relations, $\ker m\subseteq g^{-1}(\ker f)$, so
$\Hom(T_i,M)\cong\ker m\cap g^{-1}(\ker l)$. Hence
$$\frac{\Hom(T_i,M)}{\Hom(S_j,M)}\cong\frac{\ker m\cap g^{-1}(\ker
l)}{\ker m\cap\ker g}\cong g(\ker m)\cap\ker l\subseteq\ker
f\cap\ker l\cong\Hom(S_i,M),$$
and so we see that
$$\dim\Hom(S_i,M)+\dim\Hom(S_j,M)\geq\dim\Hom(T_i,M)$$
with equality if and only if the inclusion
\begin{equation*}
g(\ker m)\cap\ker l\subseteq\ker f\cap\ker l\tag{$\dagger$}
\end{equation*}
is an equality. Likewise,
$$\dim\Hom(S_i,M)+\dim\Hom(S_j,M)\geq\dim\Hom(T_j,M)$$
with equality if and only if the inclusion
\begin{equation*}
f(\ker l)\cap\ker m\subseteq\ker g\cap\ker m\tag{$\ddag$}
\end{equation*}
is an equality.

At this point, it remains to show that the inclusions ($\dagger$)
and ($\ddag$) cannot be both strict. We proceed by way of
contradiction and assume that there are
$$v\in\bigl(\ker f\cap\ker l\bigr)\setminus
\bigl(g(\ker m)\cap\ker l\bigr)\quad\text{and}\quad
w\in\bigl(\ker g\cap\ker m\bigr)\setminus
\bigl(f(\ker l)\cap\ker m\bigr).$$
We will construct a one-dimensional family $(M_t)$ of points in
$\Lambda(\nu)$ such that $M_0=M$ and
$$D_{-s_j\omega_j}(M_t)=\dim\Hom(S_j,M_t)\neq
\dim\Hom(S_j,M)=D_{-s_j\omega_j}(M)=D_{-s_j\omega_j}(Z)$$
for all $t\neq0$. Thus $M_t$ cannot be in $\Omega_Z$ for $t\neq0$,
which contradicts the fact that $M$ was chosen in the open set $\Omega_Z$.

We distinguish two cases. Suppose first that $w\notin\im f$. Then we
can choose $W_1\subseteq W$ complementary to $\spn_K(w)$ and such
that $\im f\subseteq W_1$. For each $t\in K$, we define $M_t$ to be
the same as $M$ except that the linear map $g$ is deformed to $g_t$
defined as follows:
$$g_t\bigl|_{W_1}=g\bigl|_{W_1},\quad g_t(w)=tv.$$
The preprojective relations are still satisfied, so that
$M_t\in\Lambda(\nu)$. On the other hand, for all $t\neq0$,
$$\ker g_t\cap\ker m\subsetneq\ker g\cap\ker m,\quad\text{ hence}
\quad\dim\Hom(S_j,M_t)<\dim\Hom(S_j,M).$$

Now suppose that $w\in\im f$ and choose $u\in f^{-1}(w)$. Choose a
complement $W_1$ to $\spn_K(w)$ inside $W$ such that $f(\ker l)\subseteq
W_1$. Since $w\in\ker g$, we have $gf(u)=0$ and so $kl(u)=0$. Let
$x=l(u)\in R$. Note that $x\notin l(f^{-1}(W_1))$ since $w\notin
W_1+f(\ker l)$. Hence we can find a complementary subspace $R_1$ to
$\spn_{K}(x)$ inside $R$ such that $l(f^{-1}(W_1))\subseteq R_1$.
For each $t\in K$, we define a point $M_t$ by setting
$$g_t\bigl|_{W_1}=g\bigl|_{W_1},\quad g_t(w)=tv,\quad
k_t\bigl|_{R_1}=k\bigl|_{R_1},\quad k_t(x)=\varepsilon(c)tv$$
and leaving all the other arrows unchanged. It is easy to verify the
preprojective relations for $M_t$ and to check that again, for
all $t\neq0$,
$$\ker g_t\cap\ker m\subsetneq\ker g\cap\ker m,\quad\text{ hence}
\quad\dim\Hom(S_j,M_t)<\dim\Hom(S_j,M).$$
\end{proof}

And now we complete the induction.

\trivlist
\item[\hskip\labelsep{\itshape Proof of Theorem \ref{th:TropPluRel}.}]
\upshape
Fix a pair $(i,j)$ of connected vertices.  We will show, by induction
on $\ell(w)$, that the tropical
Pl\"ucker relation holds at the hexagon $(w,i,j)$ for any $Z$.

The case $w=1$ was taken care of in Lemma~\ref{le:TropPluRelBase},
so let $w\neq1$, let $\nu$ be a dimension-vector, and let $Z$ be
an irreducible component of $\Lambda(\nu)$.

We can find $k\in I$ such that $\ell(s_kw)<\ell(w)$. We can find
$b\in B(-\infty)$ such that $Z=\Lambda_b$. Set
$$c=\varphi_k(b^*),\quad
b'=\bigl(S_k\bigl(\tilde f_k^{\max{}}b^*\bigr)\bigr){}^*\quad
\text{and}\quad Z'=\Lambda_{b'}.$$
By Proposition~\ref{pr:ReflFuncGene}, the set
$\bigl\{(M_a)\in Z\cap(\Lambda(\nu)_{k,c})^*\bigm|
[\Sigma_k^*M]\subseteq\Omega_{Z'}\bigr\}$ contains a dense open
subset of $Z$. We can thus find $M$ that belongs to both this
set and $\Omega_Z$. In particular, we have
$D_{s_k\gamma}(M)=D_{s_k\gamma}(Z)$ and
$D_\gamma(\Sigma_k^*M)=D_\gamma(Z')$
for each chamber weight $\gamma$. Combining this with
equation~(\ref{eq:InductRelPrep}), we obtain
\begin{equation} \label{eq:ZtoZ'}
D_{s_k\gamma}(Z)=D_\gamma(Z')-c\,\langle\gamma,\alpha_k\rangle
\end{equation}
whenever $\langle\gamma,\alpha_k\rangle\leq0$.

Now $(s_kw,i,j)$ is an hexagon, and the chamber weights involved
in the tropical Pl\"ucker relation for $(w,i,j)$ are of the form
$s_k\gamma$, where $\gamma$ is a weight involved in the tropical
Pl\"ucker relation for $(s_kw,i,j)$. Moreover each such weight $\gamma$
satisfies $\langle\gamma,\alpha_k\rangle\leq0$.  Hence
equation~(\ref{eq:ZtoZ'}) holds each weight $ \gamma $ involved
in the tropical Pl\"ucker relation.  From this observation, a simple 
computation shows then that the tropical Pl\"ucker relation for $Z$ at the
hexagon $(w,i,j)$ is equivalent to the tropical relation for $Z'$
at the hexagon $(s_kw,i,j)$. By the inductive hypothesis, the latter
holds, hence the former also holds. This concludes the proof.
\nobreak\noindent$\square$\endtrivlist

\subsection{MV polytopes and crystal operations}
\label{ss:MVPolyCrysOp}
Let $\MV$ be the set of all MV polytopes. Theorem~7.2 in~\cite{Kamnitzer05}
provides a natural bijection $b\mapsto P(b)$ from $B(-\infty)$ onto $\MV$.
Despite appearences, the proof does not fundamentally relies on
constructions in the affine Grassmannian; a short recapitulation of the
arguments is presented below in Remark~\ref{rk:LusParMVPol}. Under this bijection, $\MV$ inherits from $B(-\infty)$ a crystal structure.
This structure is described in section~3.6 in~\cite{Kamnitzer07}; we now
recall how it works.

Let $P\in\MV$ with vertex datum $(\mu_w)_{w\in W}$ and hyperplane
datum $(A_\gamma)_{\gamma\in\Gamma}$. Then the weight of $P$ is
$\wt(P)=\mu_{w_0}$. Fix now $i\in I$. As explained in section
\ref{ss:PWPoly}, the length of the edge between the vertices
$\mu_{w_0}$ and $\mu_{s_iw_0}$, that is, the number $c$ such that
$\mu_{w_0}-\mu_{s_iw_0}=c\alpha_i$, is the integer
$$c=A_{\omega_i}+A_{s_i\omega_i}+\sum_{\substack{j\in I\\j\neq i}}
a_{ij}A_{\omega_j}.$$
We can then define $\varphi_i(P)=c$ and
$\varepsilon_i(P)=c-\langle\alpha_i,\wt(P)\rangle$.

The key to understand the operators $\tilde e_i$ and $\tilde f_i$ is
the observation that when the $A_\gamma$ are known for all $\gamma$ in
$$\{\omega_i\mid i\in I\}\cup\{\gamma\in\Gamma\mid
\langle\gamma,\alpha_i\rangle\leq0\},$$
then the tropical Pl\"ucker relations determine the remaining $A_\gamma$.
Therefore the demand that the operators $\tilde e_i$ and $\tilde f_i$ do
not change the $A_\gamma$ such that $\langle\gamma,\alpha_i\rangle\leq0$
and the conditions
$$\wt\bigl(\tilde e_iP)=\wt(P)+\alpha_i\quad\text{and}\qquad
\wt\bigl(\tilde f_iP)=\wt(P)-\alpha_i$$
fully specify $\tilde e_iP$ and $\tilde f_iP$. Since the crystal
$\MV$ is isomorphic to $B(-\infty)$, the polytope $\tilde e_iP$ always
exists, and the polytope $\tilde f_iP$ exists if and only if
$\varphi_i(P)>0$.

By Theorem~\ref{th:TropPluRel}, to each component $Z\in\Irr\Lambda(\nu)$
is associated an MV polytope $\Pol(Z)$, defined by the BZ datum
$(D_\gamma(Z))_{\gamma\in\Gamma}$. So we have a map $\Pol$ from
$\mathbf B=\bigsqcup_{\nu\in Q_+}\Irr\Lambda(\nu)$ to~$\MV$.

\begin{theorem}
\label{th:CrysIso}
In the triangle
$$\xymatrix@R=2.4em@C=1em{&B(-\infty)\ar[dl]_(.6){b\mapsto\Lambda_b}
\ar[dr]^(.6){b\mapsto P(b)}&\\\mathbf B\ar[rr]_\Pol&&\MV,}$$
each arrow is an isomorphism of crystals. The triangle commutes.
\end{theorem}
\begin{proof}
It is enough to prove that the horizontal arrow is a morphism of
crystals, since the two diagonal arrows are isomorphism of crystals and
the identity is the only endomorphism of the crystal $B(-\infty)$.
So we are reduced to show that $\Pol$ preserves the weight map $\wt$,
the functions $\varphi_i$ and the operators $\tilde f_i^{\max{}}$.
Let $Z$ be an element of $\mathbf B$.

We have seen in section~\ref{ss:GGMSStrata} that for any $\Pi(Q)$-module
$T$, the vertex $\mu_{w_0}(T)$ of $\Pol(T)$ is equal to the
dimension-vector of $T$. Taking $T$ in $\Omega_Z$, we see that $\Pol$
preserves the weight map.

Let $i\in I$. For any $\Pi(Q)$-module $T$, we compute, using
Proposition~\ref{pr:CompRel} at the fourth step and Remark
\ref{rk:DGammaDom}~\ref{it:DGDa} at the fifth:
\begin{align*}
\dim\hd_iT&=\dim\Hom_{\Pi(Q)}(T,S_i)\\
&=\dim\Hom_{\Pi(Q)}(S_i,T^*)\\
&=D_{-s_i\omega_i}(T^*)\\
&=D_{s_i\omega_i}(T)-\langle s_i\omega_i,\dimvec T\rangle\\
&=D_{s_i\omega_i}(T)+D_{\omega_i}(T)+
\sum_{\substack{j\in I\\j\neq i}}a_{ij}D_{\omega_j}(T).
\end{align*}
Hence the function $T\mapsto\dim\hd_iT$ is constant on
$\Omega_Z$, with value
$$D_{s_i\omega_i}(Z)+D_{\omega_i}(Z)+
\sum_{\substack{j\in I\\j\neq i}}a_{ij}D_{\omega_j}(Z).$$
This number is thus the general value on $Z$ of that function, so
is equal to $\varphi_i(Z)$, by definition of the crystal structure
of $\mathbf B$. On the other hand, this expression is precisely
$\varphi_i(\Pol(Z))$, by definition of the crystal structure on
$\MV$. Therefore $\Pol$ preserves the function $\varphi_i$.

Let again $i\in I$, and set $Z'=\tilde f_i^{\max{}}Z$. By
Proposition~\ref{pr:ReformFimax}, the set $\{(M_a)\in
Z\mid[\Sigma_i^*\Sigma_iM]\subseteq\Omega_{Z'}\}$ is open and
dense in $Z$. We can thus find $M$ that belongs to both this
set and $\Omega_Z$. In particular, we have $D_\gamma(M)=D_\gamma(Z)$
and $D_\gamma(\Sigma_i^*\Sigma_iM)=D_\gamma(Z')$ for each chamber
weight $\gamma$. Equation~(\ref{eq:BZDataRadi}) then tells us that
$D_\gamma(Z)=D_\gamma(Z')$ whenever $\langle\gamma,\alpha_i\rangle\leq0$.
Since moreover
$$\wt(\Pol(Z'))=\wt(Z')=\wt(Z)-\varphi_i(Z)\alpha_i=\wt(\Pol(Z))-
\varphi_i(\Pol(Z))\alpha_i=\wt(\tilde f_i^{\max{}}\Pol(Z)),$$
we have $\Pol(Z')=\tilde f_i^{\max{}}\Pol(Z)$. Thus $\Pol$
preserves the operation $\tilde f_i^{\max{}}$.
\end{proof}

\section{From $KQ$ to $\Pi(Q)$ and back}
\label{se:FromKQToPiQ}
In this section, we investigate the relations between $KQ$-modules
and $\Pi(Q)$-modules. This theme admits several variations: we
relate our reflection functors with the classical version for
$KQ$-modules due to Bernstein, Gelfand, Ponomarev; we show that
(some of) the $\Pi(Q)$-modules $N(-w\omega_j)$ are (very close to being)
induced from an indecomposable $KQ$-module; and we determine
combinatorially the bijection from the set of $G(\nu)$-orbits in
the representation space $\Rep(KQ,\nu)$ to $\Irr\Lambda(\nu)$ that
maps an orbit $\mathscr O$ to the closure $\overline{T_\mathscr O^*}$
of its conormal bundle.

\subsection{Recall on quivers}
\label{ss:RecallQuiv}
In this section, $Q=(I,E)$ is an arbitrary quiver, not necessarily
of Dynkin type.

The lattice $\mathbb Z^I$ is then endowed with a nonsymmetric bilinear
form defined by
$$\langle\mu,\nu\rangle_Q=\sum_{i\in I}\mu_i\nu_i-
\sum_{a\in E}\mu_{s(h)}\nu_{t(h)}$$
and called the Euler form. Then the symmetric bilinear form defined
in section~\ref{ss:PreprojAlg} is
$$\langle\mu,\nu\rangle=\langle\mu,\nu\rangle_Q+\langle\nu,\mu\rangle_Q.$$
The following formula is well-known:
if $M$ and $N$ are two $KQ$-modules, then
$$\langle\dimvec M,\dimvec N\rangle=\dim\Hom_{KQ}(M,N)-
\dim\Ext^1_{KQ}(M,N).$$

Assume that $i$ is a source of $Q$ and let $\sigma_iQ$ denote the quiver
obtained from $Q$ by replacing any arrow $a\in E$ that terminates at
$i$ by the arrow $a^*$ with the opposite orientation. In this context,
we have the traditional Bernstein-Gelfand-Ponomarev reflection
functors defined in \cite{Bernstein-Gelfand-Ponomarev73}:
$$\xymatrix@C=3.5em{KQ\text{-mod}\ar@<.5ex>[r]^(.48){\Phi_i^-}&
K(\sigma_iQ)\text{-mod}.\ar@<.5ex>[l]^(.52){\Phi_i^+}}$$

In particular, $ \Phi_i^-(M) $ is defined by replacing the vector space
$ M_i $ by the cokernel of the outgoing maps from vertex $ i $ and 
$ \Phi_i^+(M) $ by replacing $ M_i $ with the kernel of the incoming maps
to vertex $ i $.  The functors $ \Phi_i^- $ and $ \Phi_i^+ $ enjoy 
properties similar to the ones stated in
Proposition~\ref{pr:ReflFuncAdj} for the functors
$\Sigma_i$ and $\Sigma_i^*$.

Certainly $KQ$ can be viewed as a subalgebra of $\Pi(Q)$. By restriction,
a $\Pi(Q)$-module $T$ gives rise to a $KQ$-module, which we denote
by $T\bigl|_Q$; in other words, we forget the action of the arrows
$a\in E^*$. Now suppose that $i$ is a source in $Q$. The quiver
$\sigma_iQ$ also gives rise to a preprojective algebra $\Pi(\sigma_iQ)$;
it has the same set $H$ of arrows as $\Pi(Q)$, but with a
different function $\varepsilon$. There is an isomorphism
$\Upsilon_i:\Pi(\sigma_iQ)\to\Pi(Q)$, which changes any arrow $b$ that
terminates at $i$ into $-b$ and which fixes all the other arrows in $H$.
We can then pull-back a $\Pi(Q)$-module $T$ by $\Upsilon_i$ and restrict
the $\Pi(\sigma_iQ)$-module $\Upsilon_i^*T$ thus obtained to $K(\sigma_iQ)$.

\begin{proposition}
\label{pr:ReflFuncBGP}
Let $T$ be a $\Pi(Q)$-module. Then
$$\Phi_i^+\Bigl((\Upsilon_i^*T)\bigl|_{\sigma_iQ}\Bigr)
\cong(\Sigma_iT)\bigl|_Q\quad\text{ and }\quad
\Phi_i^-\Bigl(T\bigl|_Q\Bigr)\cong(\Upsilon_i^*\Sigma_i^*T)
\bigl|_{\sigma_iQ}.$$
\end{proposition}
\begin{proof}
The definition of the functors $\Sigma_i$ and $\Sigma_i^*$
obviously extends Bernstein, Gelfand and Ponomarev's construction.
\end{proof}

To conclude this section, we present a lemma also related to reflection
functors for the single quiver.
\begin{lemma}
\label{le:ReflFuncEulerForm}
Assume that $i$ is a source in $Q$. Then
$\langle s_i\mu,s_i\nu\rangle_Q=\langle\mu,\nu\rangle_{\sigma_iQ}$
for all $\mu,\nu\in\mathbb Z^I$.
\end{lemma}
\begin{proof}
Using the definition of the Euler form and that $i$ is a source in $Q$,
one finds that
$$\langle\mu-\mu_i\alpha_i,\nu\rangle_Q=
\sum_{\substack{j\in I\\j\neq i}}\mu_j\nu_j-\sum_{\substack{a\in E\\
i\notin\{s(a),t(a)\}}}\mu_{s(a)}\nu_{t(a)}=
\langle\mu,\nu-\nu_i\alpha_i\rangle_{\sigma_iQ}.$$
One then computes
\begin{align*}
\langle s_i\mu,\nu\rangle_Q
&=\langle\mu,\nu\rangle_Q-\langle\alpha_i,\mu\rangle\,\langle\alpha_i,
\nu\rangle_Q\\
&=\langle\mu,\nu\rangle_Q-\bigl(\langle\mu,\alpha_i\rangle_Q
+\langle\alpha_i,\mu\rangle_Q\bigr)\,\langle\alpha_i,\nu\rangle_Q\\
&=\langle\mu-\mu_i\alpha_i,\nu\rangle_Q-\langle\alpha_i,\mu\rangle_Q
\;\langle\alpha_i,\nu\rangle_Q\\
&=\langle\mu,\nu-\nu_i\alpha_i\rangle_{\sigma_iQ}-\langle\mu,
\alpha_i\rangle_{\sigma_iQ}\;\langle\nu,\alpha_i\rangle_{\sigma_iQ}\\
&=\langle\mu,\nu\rangle_{\sigma_iQ}-\bigl(\langle\alpha_i,\nu
\rangle_{\sigma_iQ}+\langle\nu,\alpha_i\rangle_{\sigma_iQ}\bigr)
\;\langle\mu,\alpha_i\rangle_{\sigma_iQ}\\
&=\langle\mu,\nu\rangle_{\sigma_iQ}-\langle\alpha_i,\nu\rangle\,
\langle\mu,\alpha_i\rangle_{\sigma_iQ}\\
&=\langle\mu,s_i\nu\rangle_{\sigma_iQ}.
\end{align*}
This formula is of course equivalent to the one written in the statement.
\end{proof}

\begin{other*}{Remark}
There is a more conceptual proof of this lemma. The
dimension-vector defines an isomorphism between $\mathbb Z^I$ and
the Grothendieck group of the derived category
$\mathbb D^b(KQ\text{-mod})$, and one has
$$\sum_{n\in\mathbb Z}(-1)^n\dim\Ext^n_{\mathbb D^b(K(Q))}(M,N)
=\langle\dimvec M,\dimvec N\rangle_Q$$
for any two objects $M$ and $N$ in this category.
A similar statement applies to the derived category
$\mathbb D^b(K(\sigma_iQ)\text{-mod})$. Now the right-exact
functor $\Phi_i^+$ induces an equivalence of derived categories
$$\mathbb R\Phi_i^+:\mathbb D^b(K(\sigma_iQ)\text{-mod})\to\mathbb
D^b(KQ\text{-mod})$$
and one has
$$\dimvec\mathbb R\Phi_i^+M=s_i\left(\dimvec M\right).$$
\end{other*}

\subsection{Modules $N(\gamma)$ and induction}
\label{ss:ModNGammaInd}
We now go back to our set-up where $Q$ is obtained by choosing
an orientation of the Dynkin diagram~$\Gamma$.

We denote the quiver opposite to $Q$ by $Q^*$. Recall that a sequence
$(i_1,\ldots,i_n)$ of vertices of is called adapted to $Q^*$ if $i_1$
is a source of $Q$, $i_2$ is a source of $\sigma_{i_1}Q$, $i_3$ is a
source of $\sigma_{i_2}\sigma_{i_1}Q$, and so on (see section~4.7 in
\cite{Lusztig90a}).

We first prove a combinatorial version of our induction result.

\begin{proposition}
\label{pr:CombiInduction}
Assume that the sequence $(i_1,\ldots,i_n)$ is adapted to $Q^*$ and that
$s_{i_1}\cdots s_{i_n}$ is a reduced expression of an element $w\in W$.
Then $\langle-w\alpha_{i_n},?\rangle_Q=\langle-w\omega_{i_n},?\rangle$ as
linear forms on $Q$.
\end{proposition}
\begin{proof}
For each $j\in I$ different from $i$,
$$\langle\alpha_{i_1}-\omega_{i_1},\alpha_j\rangle=
\langle\alpha_{i_1},\alpha_j\rangle=
\langle\alpha_{i_1},\alpha_j\rangle_Q+\langle\alpha_{i_1},
\alpha_j\rangle_Q=\langle\alpha_{i_1},\alpha_j\rangle_Q,$$
because $i_1$ is a source in $Q$. The equality
$\langle\alpha_{i_1}-\omega_{i_1},\alpha_j\rangle=
\langle\alpha_{i_1},\alpha_j\rangle_Q$
also holds for $j=i_1$. We deduce by linearity that
$$\langle\alpha_{i_1}-\omega_{i_1},?\rangle=
\langle\alpha_{i_1},?\rangle_Q,$$
which is the case $n=1$.

For higher value of $n$, we proceed by induction, using Lemma
\ref{le:ReflFuncEulerForm} and the $W$-invariance of $\langle\;,\;\rangle$.
\end{proof}

We will use the notation $S_i$ for the simple $KQ$-module attached to the
vertex $i$ and we will also speak of $i$-socle
for $KQ$-modules.

Recall that by the classical theorem of Gabriel, for each positive root 
$\alpha$, there is a unique (up to isomorphism) indecomposable
$KQ$-module with dimension-vector $ \alpha$, which we denote by 
$M(\beta)$.  Moreover, 
Berstein-Gelfand-Ponomarev \cite{Bernstein-Gelfand-Ponomarev73} showed
that these indecomposable modules may be constructed using their
reflection functors.  Let $ Q $ be a quiver and let $i $ be a source in
$Q$.  Let $ \beta $ be a positive weight.  Then
\begin{equation} \label{eq:reflectindecomp}
\Phi_i^- M(\beta) = \begin{cases}
0 \text{ if } \beta = \alpha_i \\
M(s_i \beta) \text{ otherwise}
\end{cases}
\end{equation}
where $M(s_i \beta) $ denotes the indecomposable $ K(\sigma_i Q) $ module
with dimension-vector $ s_i \beta $.  Starting with the base case
$M(\alpha_i) = S_i $, this gives a recursive way to construct all
$M(\beta) $.

Ringel
has observed (see section~4 in \cite{Ringel96}) that if $\alpha$ and
$\beta$ are two positive roots, then at least one of the two spaces
$\Hom_{KQ}(M(\alpha),M(\beta))$ or $\Ext^1_{KQ}(M(\alpha),M(\beta))$
is zero. With this information in hand, the module-theoretic expression
of the Euler form writes
\begin{equation*}
\left\{\begin{aligned}
\dim\Hom_{KQ}(M(\alpha),M(\beta))
=\max(0,\langle\alpha,\beta\rangle_Q)\\[4pt]
\dim\Ext^1_{KQ}(M(\alpha),M(\beta))
=\max(0,-\langle\alpha,\beta\rangle_Q).
\end{aligned}\right.
\end{equation*}

The following proposition is the analogue for the single quiver $Q$ of
equation~(\ref{eq:InductRelPrep}).
\begin{proposition}
\label{pr:InductRelSing}
Assume that $i$ is a source in $Q$ and let $\beta$ be a positive root,
different from $\alpha_i$. Let $T$ be a $KQ$-module. Then
$$\dim\Hom_{KQ}(M(s_i\beta),T)=\dim\Hom_{K(\sigma_iQ)}(M(\beta),
\Phi_i^-T)-\langle\beta,\dimvec\soc_iT\rangle_{\sigma_iQ}.$$
(In this equality, it is implicitly understood that $M(s_i\beta)$ is
the indecomposable $KQ$-module with dimension-vector $s_i\beta$, and
that $M(\beta)$ is the indecomposable $K(\sigma_iQ)$-module with
dimension-vector $\beta$.)
\end{proposition}
\begin{proof}
We first notice that $\langle s_i\beta,\alpha_i\rangle_Q\geq0$, for
$i$ is a source in $Q$ and $s_i\beta\in Q_+$. Thus
$$\Ext^1_{KQ}(M(s_i\beta),S_i)=0\quad\text{ and }\quad
\dim\Hom_{KQ}(M(s_i\beta),S_i)=\langle s_i\beta,\alpha_i\rangle_Q.$$
Lemma~\ref{le:ReflFuncEulerForm} moreover says that $\langle
s_i\beta,\alpha_i\rangle_Q=-\langle\beta,\alpha_i\rangle_{\sigma_iQ}$.

On the other hand, the functor $\Phi_i^-$ maps $M(s_i\beta)$ to
$M(\beta)$ by (\ref{eq:reflectindecomp}); by adjunction, we thus get
$$\Hom_{K(\sigma_iQ)}(M(\beta),\Phi_i^-T)\cong
\Hom_{K(\sigma_iQ)}(\Phi_i^-M(s_i\beta),\Phi_i^-T))\cong
\Hom_{KQ}(M(s_i\beta),\Phi_i^+\Phi_i^-T).$$
It remains to apply the functor $\Hom_{KQ}(M(s_i\beta),?)$ to the
(split) short exact sequence
$$0\to\soc_iT\to T\to\Phi_i^+\Phi_i^-T\to0$$
and to take dimensions.
\end{proof}

We now have all the ingredients to prove the following result.
\begin{theorem}
\label{th:ModInduction}
Assume that the sequence $(i_1,\ldots,i_n)$ is adapted to $Q^*$ and that
$s_{i_1}\cdots s_{i_n}$ is a reduced expression of an element $w\in W$.
Then for any $\Pi(Q)$-module $T$,
$$\dim\Hom_{KQ}(M(-w\alpha_{i_n}),T\bigl|_Q)=
\dim\Hom_{\Pi(Q)}(N(-w\omega_{i_n}),T).$$
\end{theorem}
\begin{proof}
We first observe that since $i_1$ is a source of $Q$, the
$i_1$-socles of the $\Pi(Q)$-module $T$ and of the $KQ$-module
$T\bigl|_Q$ coincide. In particular, they have the same
dimension (-vector).

The case $n=1$ is then obvious, because then both $M(-w\alpha_{i_n})$
and $N(-w\omega_{i_n})$ are the simple module $S_{i_1}$.

In the case $n>1$, we proceed by induction on $n$. More precisely, we
assume that the equality holds for the sequence $(i_2,\ldots,i_n)$,
the quiver $\sigma_{i_1}Q$ and the $\Pi(\sigma_{i_1}Q)$-module
$\Upsilon_{i_1}^*\Sigma_{i_1}^*T$. A straightforward computation
based on Propositions~\ref{pr:BZDataPrep}, \ref{pr:ReflFuncBGP},
\ref{pr:CombiInduction} and~\ref{pr:InductRelSing} then leads to
the equality for the sequence $(i_1,\ldots,i_n)$, the quiver~$Q$ and
the $\Pi(Q)$-module~$T$.
\end{proof}

Theorem~\ref{th:ModInduction} shows the existence of an isomorphism
$$\Hom_{KQ}(M(-w\alpha_{i_n}),T\bigl|_Q)\cong
\Hom_{\Pi(Q)}(N(-w\omega_{i_n}),T).$$
We believe that this isomorphism can be made functorial in
$T\in\Pi(Q)$-mod; in other words, we believe that $N(-w\omega_{i_n})$
is isomorphic to $\Pi(Q)\otimes_{KQ}M(-w\alpha_{i_n})$. This fact
can probably be deduced from the work of Gei\ss, Leclerc and Schr\"oer
\cite{Geiss-Leclerc-Schroer10}. Anyway, this would generally not give
a description of all the modules $N(\gamma)$ with $\gamma$ a chamber
weight; indeed only in type $A$ can one always write
$\gamma=-s_{i_1}\cdots s_{i_n}\omega_{i_n}$ for a sequence
$(i_1,\ldots,i_n)$ adapted with an orientation of the graph $\Gamma$.

\subsection{Irreducible components and conormal bundles}
\label{ss:IrrCompConoBund}
Let $\nu\in Q_+$ be a dimension-vector. The isomorphism classes of
representations of the
quiver $Q$ with dimension-vector $\nu$ are the orbits of the
the group $G(\nu)$ on the representation space $\Rep(KQ,\nu)$. In
Proposition~14.2 of~\cite{Lusztig91}, Lusztig shows that the
irreducible components of $\Lambda(\nu)$ are precisely the closures
of the conormal bundles to the $G(\nu)$-orbits in $\Rep(KQ,\nu)$.
More precisely, he observes (\textit{loc.~cit.}, 12.8~(a)) that
under the trace duality
$$\prod_{a\in H}\Hom_K(K^{\nu_{s(a)}},K^{\nu_{t(a)}})\cong
\Rep(KQ,\nu)\times\Rep(KQ,\nu)^*,$$
any point $(M_a)\in\Lambda(\nu)$ belongs to the conormal bundle
$T^*\mathscr O$, where $\mathscr O$ is the $G(\nu)$-orbit in
$\Rep(KQ,\nu)$ that corresponds to the isomorphism class of the
$KQ$-module $M\bigl|_Q$.

Lusztig's interest in the varieties $\Lambda(\nu)$ arose in connection
with its study of the canonical basis $\mathcal B$ of the quantum
group $U_q(\mathfrak n)$, the positive part of $U_q(\mathfrak g)$.
Geometrically, an element $b$ in the canonical basis is represented
by a simple $G(\nu)$-equivariant perverse sheaf $L_{b,Q}$ on
$\Rep(KQ,\nu)$, where $\nu$ is the weight of $b$, and such a sheaf is
always the intersection cohomology sheaf
$\IC\bigl(\overline{\mathscr O_{b,Q}}\bigr)$ (suitably shifted)
of the closure of a $G(\nu)$-orbit (with coefficients in the trivial
local system). The singular support $SS(L_{b,Q})$ then contains
$\overline{T^*\mathscr O_{b,Q}}$.

The canonical basis $\mathcal B$ has the structure of a crystal,
which identifies it with $B(-\infty)$. We can thus change slightly
the notation and suppose that the letter $b$ in $L_{b,Q}$ and
$\mathscr O_{b,Q}$ denotes an element of $B(-\infty)$.  It is natural
to guess the labelling of the canonical and semicanonical bases 
are compatible in the sense that 
$ \Lambda_b = \overline{T^*\mathscr O_{b,Q}}$.  This would imply the
inclusion $\Lambda_b\subseteq SS(L_{b,Q})$.  
In their article \cite{Kashiwara-Saito97}, Kashiwara and Saito
attribute this inclusion to
\cite{Lusztig90a}, but we were not able to locate this statement (or 
its proof) in Lusztig's paper.

Our aim in this final section is to apply our results to prove the 
equality $ \Lambda_b = \overline{T^*\mathscr O_{b,Q}}$.
In other words, we want to prove that the map
$b\mapsto\overline{T^*\mathscr O_{b,Q}}$ from $B(-\infty)$ to
$\mathbf B$ is a morphism of crystal. The main difficulty is
that one does not a priori know that
$\overline{T^*\mathscr O_{b,Q}}=\overline{T^*\mathscr O_{b,Q'}}$
for two different orientations $Q$ and~$Q'$.

We first need to recall some facts concerning parametrizations of
$B(-\infty)$.  As before, $w_0$ is the longest
element in $W$. We denote its length by $N$ and call $\mathscr X$
the set of all tuples $\mathbf i=(i_1,\ldots,i_N)$ such that
$s_{i_1}\cdots s_{i_N}$ is a reduced expression of $w_0$. As
explained by Lusztig \cite{Lusztig90a}, to each $\mathbf i\in\mathscr X$
corresponds a bijection $\mathbf n(?,\mathbf i):B(-\infty)\to\mathbb N^N$,
which is the combinatorial counterpart of the transition matrix between
the canonical basis $\mathcal B$ and the PBW basis of $U_q(\mathfrak n)$
defined by $\mathbf i$. These bijections can be collectively
characterized by the properties (L1), (L2a) and (L2b) below.
In these statements, $(b,\mathbf i)\in B(-\infty)\times\mathscr X$
and $(n_1,\ldots,n_N)$ is the tuple $\mathbf n(b,\mathbf i)$.
\begin{description}
\item[(L1)]
One has
$$n_1=\varphi_{i_1}(b),\quad
\mathbf n(\tilde e_{i_1}b,\mathbf i)=(n_1+1,n_2,\ldots,n_N),\quad
\mathbf n(\tilde f_{i_1}^{\max{}}b,\mathbf i)=(0,n_2,\ldots,n_N).$$
\item[(L2a)]
Let $(i,j)\in I^2$ and $k\in\{1,\ldots,N-1\}$ be such that $a_{ij}=0$
and $(i_k,i_{k+1})=(i,j)$. Set
$\mathbf j=(i_1,\ldots,i_{k-1},j,i,i_{k+2},\ldots,i_N)$. Then
$$\mathbf n(b,\mathbf j)=
(n_1,\ldots,n_{k-1},n_{k+1},n_k,n_{k+2},\ldots,n_N).$$
\item[(L2b)]
Let $(i,j)\in I^2$ and $k\in\{2,\ldots,N-1\}$ be such that $a_{ij}=-1$
and $(i_{k-1},i_k,i_{k+1})=(i,j,i)$. Let $(p,q,r)=(n_{k-1},n_k,n_{k+1})$
and set $\mathbf j=(i_1,\ldots,i_{k-2},j,i,j,i_{k+2},\ldots,i_N)$,
$q'=\min(p,r)$ and $(p',r')=(r+q-q',p+q-q')$. Then
$$\mathbf n(b,\mathbf j)=
(n_1,\ldots,n_{k-2},p',q',r',n_{k+2},\ldots,n_N).$$
\item[(L3)]
Let $i_1^*\in I$ be the element such that $w_0\alpha_{i_1}=-\alpha_{i_1^*}$.
Suppose $n_1=0$ and set $\mathbf j=(i_2,i_3,\ldots,i_N,i_1^*)$. Then
$$\mathbf n(S_{i_1}(b),\mathbf j)=(n_2,n_3,\ldots,n_N,0).$$
\end{description}

Property (L3) is due to Saito (Proposition~3.4.7 in \cite{Saito94}).
The interesting point here is that with the help of (L1) and (L3),
one can determine $\mathbf n(b,\mathbf i)$ by applying to $b$ the
operations $\varphi_i$, $\tilde f_i^{\max{}}$ and $S_i$, with $i$
running over the successive terms of the sequence $\mathbf i$. This
procedure is convenient for us, for these three crystal operations
have a clear meaning in the preprojective model, thanks to
Propositions~\ref{pr:ReformFimax} and~\ref{pr:ReflFuncGene}.

\begin{other}{Remark}
\label{rk:LusParMVPol}
In this remark, we explain the origin of the bijection $B(-\infty)\cong\MV$
we used in section~\ref{ss:MVPolyCrysOp}.
\begin{enumerate}
\item\label{it:LPMPa}
For $\mathbf i,\mathbf j\in\mathscr X$, we can define a permutation
$R_{\mathbf i}^{\mathbf j}$ of $\mathbb N^N$ by the condition
$R_{\mathbf i}^{\mathbf j}\circ\mathbf n(?,\mathbf i)=\mathbf
n(?,\mathbf j)$. Then $b\mapsto\mathbf n(b,?)$ is a bijection from
$B(-\infty)$ onto
$$\widehat{\mathscr X}=\bigl\{\sigma:\mathscr X\to\mathbb N^N\bigm|
\forall(\mathbf i,\mathbf j)\in\mathscr X^2,\ \sigma(\mathbf j)=
R_{\mathbf i}^{\mathbf j}(\sigma(\mathbf i))\bigr\}.$$
The permutations $R_{\mathbf i}^{\mathbf j}$ extends in a natural
fashion to $\mathbb Z^N$, for it is given by a piecewise linear
formula which also makes sense for signed integers. Let $\widetilde
{\mathscr X}$ be the set defined similarly as $\widehat{\mathscr X}$,
but with $\mathbb Z^N$ instead of $\mathbb N^N$.
\item\label{it:LPMPb}
To $\sigma\in\widetilde{\mathscr X}$, we associate a vertex datum
$(\mu_w)_{w\in W}$ as follows: for each $w\in W$, we can find
$\mathbf i\in\mathscr X$ and $k\in\{0,\ldots,N\}$ such that
$ww_0=s_{i_1}\cdots s_{i_k}$, and we set
$$\mu_w=\sum_{r=k+1}^Nn_r\beta_r,\quad\text{where}\quad
\sigma(\mathbf i)=(n_1,\ldots,n_N)\quad\text{and}\quad
\beta_r=s_{i_1}\cdots s_{i_{r-1}}\alpha_{i_r}.$$
(This weight $\mu_w$ depends solely on $\sigma$ and $w$ and not on the
choice of $\mathbf i$.) We then associate an hyperplane datum
$(A_\gamma)_{\gamma\in\Gamma}$ by the equation (\ref{eq:RelVertHyp}).
Then the map $\sigma\mapsto(A_\gamma)_{\gamma\in\Gamma}$ is a bijection
from $\widetilde{\mathscr X}$ onto the set of all families that satisfies
the conditions (BZ1) and (BZ3) from section~\ref{ss:TropPlucRel}. This
fact is proved by specializing Theorem~4.3
in~\cite{Berenstein-Zelevinsky97} to the tropical semifield.
\item\label{it:LPMPc}
The numbers that appear in the left hand side of the edge inequalities
(\ref{eq:EdgeIneq}) are precisely the components of the tuples
$\sigma(\mathbf i)$. Thus $(A_\gamma)_{\gamma\in\Gamma}$ satisfies (BZ2)
if and only if each $\sigma(\mathbf i)$ belongs to $\mathbb N^N$. By
restriction, we therefore get a bijection from $\widehat{\mathscr X}$
onto the set of all BZ data. In the end, we obtain a bijection
$B(-\infty)\cong\MV$. To sum up, the MV polytope $P(b)$ associated to
an element $b\in B(-\infty)$ packs together in a geometrical way the
numerical data $\mathbf n(b,\mathbf i)$ for all possible $\mathbf i$.
\item\label{it:LPMPd}
Keeping the same notation, we now turn to the crystal structure defined
on $\MV$ by transport from $B(-\infty)$. For each $w\in W$ such that
$\ell(s_iww_0)<\ell(ww_0)$, we can choose an $\mathbf i\in\mathscr X$
that starts with $i$ in the construction of \ref{it:LPMPb} above and
we have $k>0$; thus $n_1$ does not enter in the equation that expresses
$\mu_w$, and Property (L1) tells us that the crystal operations
$\tilde e_i$ and $\tilde f_i$ do not change $\mu_w$. In terms of the
BZ data $(A_\gamma)_{\gamma\in\Gamma}$, this translates into the fact
that the operators $\tilde e_i$ and $\tilde f_i$ do not change the
$A_\gamma$ such that $\langle\gamma,\alpha_i\rangle\leq0$. We thus
recover the criterion we used in section~\ref{ss:MVPolyCrysOp}.
\item\label{it:LPMPe}
We can understand Property (L3) from the point of view of 
Theorem~\ref{th:CrysIso}.  In more detail,
equations~(\ref{eq:BZDataRefl}) and~(\ref{eq:RelVertHyp}) together
imply the relation
$$\mu_w(T)=s_i\mu_{s_iw}(\Sigma_iT)$$
for each $\Pi(Q)$-module $T$ and each $w\in W$ such that
$\ell(s_iw)>\ell(w)$. Then (L3) is a direct translation of
Proposition~\ref{pr:ReflFuncGene}.
\end{enumerate}
\end{other}

We are now ready to prove our result.
\begin{proposition}
\label{pr:IrrCompConoBund}
The equality $\Lambda_b=T^*\mathscr O_{b,Q}$ holds for each
$b\in B(-\infty)$.
\end{proposition}
\begin{proof}
Let $\mathbf i\in\mathscr X$ be adapted to the quiver $Q$:
$i_1$ is a sink of $Q$, $i_2$ is a sink of $\sigma_{i_1}^{-1}Q$,
and so on. (Such a sequence $\mathbf i$ always exists, thanks to
Proposition~4.12~(b) in \cite{Lusztig90a}).

It is enough to show that two elements $b,b'\in B(-\infty)$ such
that $\Lambda_b=\overline{T^*\mathscr O_{b',Q}}$ have necessarily
the same image under the map $\mathbf n(?,\mathbf i)$. So let
us take two such elements and set
$$(n_1,\ldots,n_N)=\mathbf n(b,\mathbf i)\quad\text{and}\quad
(n'_1,\ldots,n'_N)=\mathbf n(b',\mathbf i).$$

Let $T$ be a general point in $\Lambda_b$. In view of Property (L3)
above, successive applications of Proposition~\ref{pr:ReflFuncGene}
show that
$$n_k=\dim\hd_{i_k}\bigl(\Sigma_{i_{k-1}}\cdots\Sigma_{i_1}T\bigr).$$
We can here substitute
$$\bigl(\Sigma_{i_{k-1}}(\Upsilon_{i_{k-1}}^{-1})^*\bigr)
\cdots\bigl(\Sigma_{i_1}(\Upsilon_{i_1}^{-1})^*\bigr)\cong
((\Upsilon_{i_{k-1}}\cdots\Upsilon_{i_1})^{-1})^*
(\Sigma_{i_{k-1}}\cdots\Sigma_{i_1})$$
to $\Sigma_{i_{k-1}}\cdots\Sigma_{i_1}$
without changing the $i_k$-head.

On the other hand, for each positive root $\beta$, let $M(\beta)$ be
the indecomposable $KQ$-module with dimension-vector $\beta$, and let
$\beta_1$, \dots, $\beta_N$ be the indexation of the positive roots
defined by $\mathbf i$. By sections~4.15--4.16 in \cite{Lusztig90a},
the fact that $T\bigl|_Q\in\mathscr O_{b',Q}$ means
$$T\cong M(\beta_1)^{\oplus n'_1}\oplus\cdots\oplus M(\beta_N)^{\oplus
n'_N}.$$
It follows from equation~(\ref{eq:reflectindecomp}) that
$$n'_k=\dim\hd_{i_k}\Bigl(\Phi_{i_{k-1}}^+\cdots\Phi_{i_1}^+
\Bigl(T\bigl|_Q\Bigr)\Bigr).$$

An appeal to Proposition~\ref{pr:ReflFuncBGP} now concludes the proof;
one has just to observe that since $i_k$ is a sink of the quiver
$Q'=\sigma_{i_{k-1}}^{-1}\cdots\sigma_{i_1}^{-1}Q$, the $i_k$-head of
the $\Pi(Q')$-module
$$\bigl(\Sigma_{i_{k-1}}(\Upsilon_{i_{k-1}}^{-1})^*\bigr)
\cdots\bigl(\Sigma_{i_1}(\Upsilon_{i_1}^{-1})^*\bigr)\bigl(T\bigr)$$
is the same as its head as a $KQ'$-module.
\end{proof}


\begin{thebibliography}{99}
\bibitem{Anderson03}
J.~E.~Anderson, A polytope calculus for semisimple groups,
\textit{Duke Math.\ J.\/} \textbf{116} (2003), 567-588.
\bibitem{Auslander-Reiten-Smalo95}
M.~Auslander, I.~Reiten, S.~O.~Smal\o, \textit{Representation theory
of Artin algebras}, Cambridge Studies in Advanced Mathematics vol.~36,
Cambridge University Press, 1995.
\bibitem{Berenstein-Zelevinsky97}
A.~Berenstein, A.~Zelevinsky, Total positivity in Schubert varities,
\textit{Comment.\ Math.\ Helv.} \textbf{72} (1997), 128--166.
\bibitem{Bernstein-Gelfand-Ponomarev73}
I.~N.~Bernstein, I.~M.~Gelfand, V.~A.~Ponomarev, Coxeter functors and
Gabriel's theorem, \textit{Uspehi Mat.\ Nauk} \textbf{28} (1973), 19--33.
\bibitem{Crawley-Boevey00}
W.~Crawley-Boevey, On the exceptional fibres of Kleinian singularities,
\textit{Amer.\ J.\ Math.\/} \textbf{122} (2000), 1027--1037.
\bibitem{Crawley-Boevey-Holland98}
W.~Crawley-Boevey, M.~P.~Holland, Noncommutative deformations of Kleinian
singularities, \textit{Duke Math.\ J.\/} \textbf{92} (1998), 605--635.
\bibitem{Geiss-Leclerc-Schroer08}
C.~Gei\ss, B.~Leclerc, J.~Schr\"oer, Partial flag varieties and preprojective
algebras, \textit{Ann.\ Inst.\ Fourier (Grenoble)\/} \textbf{58} (2008),
825--876.
\bibitem{Geiss-Leclerc-Schroer10}
C.~Gei\ss, B.~Leclerc, J.~Schr\"oer, Kac-Moody groups and cluster
algebras, arXiv:1001.3545.
\bibitem{Hiriart-Urruty-Lemarechal01}
J.-B.~Hiriart-Urruty, C.~Lemar\'echal, \textit{Fundamentals of convex
analysis}, Grundlehren Text Editions, Springer-Verlag, 2001.
\bibitem{Kamnitzer05}
J.~Kamnitzer, Mirkovic-Vilonen cycles and polytopes, \textit{Ann.\ of Math.\/}
\textbf{171} (2010), No. 1, 245--294; arXiv:math.AG/0501365.
\bibitem{Kamnitzer07}
J.~Kamnitzer, The crystal structure on the set of Mirkovi\'c-Vilonen
polytopes, \textit{Adv.\ Math.\/} \textbf{215} (2007), 66--93.
\bibitem{Kamnitzer-Sadanand10}
J.~Kamnitzer, C.~Sadanand, Modules with 1-dimensional socle and components
of Lusztig quiver varieties in type A; arXiv:1009.0272.
\bibitem{Kashiwara95}
M.~Kashiwara, On crystal bases, in: \textit{Representations of groups
(Banff, AB, 1994)}, pp.~155--197, CMS Conf.\ Proc.\ vol.~16, American
Mathematical Society, 1995.
\bibitem{Kashiwara-Saito97}
M.~Kashiwara, Y.~Saito, Geometric construction of crystal bases,
\textit{Duke Math.\ J.\/} \textbf{89} (1997), 9--36.
\bibitem{Lusztig90a}
G.~Lusztig, Canonical bases arising from quantized enveloping algebras,
\textit{J.\ Amer.\ Math.\ Soc.\/} \textbf3 (1990), 447--498.
\bibitem{Lusztig90b}
G.~Lusztig, Canonical bases arising from quantized enveloping algebras II,
\textit{Progr.\ Theoret.\ Phys.\ Suppl.} \textbf{102} (1990), 175--201.
\bibitem{Lusztig91}
G.~Lusztig, Quivers, perverse sheaves, and quantized enveloping algebras,
\textit{J.\ Amer.\ Math.\ Soc.\/} \textbf4 (1991), 365--421.
\bibitem{Lusztig96}
G.~Lusztig, Braid group action and canonical bases, \textit{Adv.\
Math.\/} \textbf{122} (1996), 237--261.
\bibitem{Lusztig98}
G.~Lusztig, On quiver varieties, \textit{Adv.\ Math.\/} \textbf{136}
(1998), 141--182.
\bibitem{Lusztig00}
G.~Lusztig, Semicanonical bases arising from enveloping algebras,
\textit{Adv.\ Math.\/} \textbf{151} (2000), 129--139.
\bibitem{Nakajima94}
H.~Nakajima, Instantons on ALE spaces, quiver varieties, and Kac-Moody
algebras, \textit{Duke Math.\ J.\/} \textbf{76} (1994), 365--416.
\bibitem{Nakajima98}
H.~Nakajima, Quiver varieties and Kac-Moody algebras,
\textit{Duke Math.\ J.\/} \textbf{91} (1998), 515--560.
\bibitem{Riedtmann80}
C.~Riedtmann, Algebren, Darstellungsk\"ocher, \"Uberlagerungen und zur\"uck,
\textit{Comment.\ Math.\ Helv.\/} \textbf{55} (1980), 199--224.
\bibitem{Ringel96}
C.~M.~Ringel, PBW-bases of quantum groups, \textit{J.\ Reine Angew.\
Math.\/} \textbf{470} (1996), 51--88.
\bibitem{Ringel98}
C.~M.~Ringel, The preprojective algebra of a quiver, in: \textit{Algebras
and modules, II (Geiranger, 1996)}, pp.~467--480, CMS Conf.\ Proc.\
vol.~24, American Mathematical Society, 1998.
\bibitem{Saito94}
Y.~Saito, PBW basis of quantized universal enveloping algebras,
\textit{Publ.\ Res.\ Inst.\ Math.\ Sci.\/} \textbf{30} (1994), 209--232.
\bibitem{Saito02}
Y.~Saito, Crystal bases and quiver varieties, \textit{Math.\ Ann.\/}
\textbf{324} (2002), 675--688.
\bibitem{Savage07}
A.~Savage, Geometric and combinatorial realizations of crystals
of enveloping algebras, in: \textit{Lie algebras, vertex operator
algebras and their applications (Raleigh, NC, 2005)}, pp.~221--232,
Contemp.\ Math.\ vol.~442, American Mathematical Society, 2007.
\end{thebibliography}
\end{document}